\documentclass[final,5p,times,twocolumn,authoryear]{elsarticle}

\usepackage{amssymb}
\usepackage{latexsym}
\usepackage{amsmath,amsthm,amsfonts}
\usepackage{diagbox}
\usepackage{algorithmic}
\usepackage{xcolor}
\usepackage[ruled,norelsize,vlined,linesnumbered]{algorithm2e}
\usepackage[caption=false,font=normalsize,labelfon
t=sf,textfont=sf]{subfig}
\usepackage{float} 

\newcommand{\rv}[1]{\textcolor{black}{#1}}
\journal{.}

\begin{document}
\onecolumn

\begin{frontmatter}



\title{A New Parallel Cooperative Landscape Smoothing Algorithm and Its Applications on TSP and UBQP}


\author[1]{Wei Wang}
\ead{ww123@stu.xjtu.edu.cn}

\author[1]{Jialong Shi\corref{cor1}}
\ead{jialong.shi@xjtu.edu.cn}

\author[1]{Jianyong Sun}
\ead{jy.sun@xjtu.edu.cn}

\author[2]{Arnaud Liefooghe}
\ead{arnaud.liefoogheuniv-littoral.fr}

\author[3]{Qingfu Zhang}
\ead{qingfu.zhang@cityu.edu.hk}

\address[1]{School of Mathematics and Statistics, Xi' an Jiaotong University, Xi' an, China}
\address[2]{LISIC, Université du Littoral Côte d'Opale, F-62228 Calais, France}
\address[3]{Department of Computer Science, City University of Hong Kong, Hong Kong, China}

\cortext[cor1]{Corresponding author}

\begin{abstract}
Combinatorial optimization problem (COP) is difficult to solve because of the massive number of local optimal solutions in his solution space. Various methods have been put forward to smooth the solution space of COPs, including homotopic convex (HC) transformation for the traveling salesman problem (TSP). \rv{This paper extends the HC transformation approach to unconstrained binary quadratic programming (UBQP) by proposing a method to construct a unimodal toy UBQP of any size. We theoretically prove the unimodality of the constructed toy UBQP. After that, we apply this unimodal toy UBQP to smooth the original UBQP by using the HC transformation framework and empirically verify the smoothing effects.} Subsequently, we introduce an iterative algorithmic framework incorporating HC transformation, referred as landscape smoothing iterated local search (LSILS). Our experimental analyses, conducted on various UBQP instances show the effectiveness of LSILS. Furthermore, this paper proposes a parallel cooperative variant of LSILS, denoted as PC-LSILS and apply it to both the UBQP and the TSP. Our experimental findings highlight that PC-LSILS improves the smoothing performance of the HC transformation, and further improves the overall performance of the algorithm.
\end{abstract}


%\item A landscape smoothing method called HC transformation is proposed for UBQP.
%\item A parallel cooperative framework is proposed for HC transformation-based algorithms.
%\item Experiments show HC transformation and parallel cooperative framework are effective.


\end{frontmatter}

\section{Introduction}
COPs are a class of problems mainly to find an optimal combination to maximize or minimize some performance metrics with limited resources or subject to some constraints. COPs are typically categorized as NP-hard, for which classical optimization methods struggle to find the optimal solution within a reasonable amount of time and become less applicable. Consequently, researchers usually use heuristics or metaheuristics to find near-optimal solutions within a reasonable amount of time. One of the key challenges associated with solving COPs is the presence of numerous local optima in the solution space, primarily attributable to the rugged and irregular nature of their fitness landscapes. It can be hypothesized that smoothing the landscapes of COPs could significantly facilitate the attainment of the global optima. In this study, we focus on two widely studied combinatorial optimization problems, namely the UBQP and the TSP, to investigate this conjecture further.

In this article, we focus on a recently proposed technique for landscape smoothing known as HC transformation, initially introduced by \cite{shi2020homotopic}. The HC transformation method, originally designed for the TSP, involves defining the HC transformation of a TSP instance as a convex combination of the original TSP and a well-designed toy TSP with a unimodal smooth landscape. The unimodal toy TSP is constructed based on known high-quality local optima, thus maintaining the information in elite solutions during the search, as illustrated in \figurename~ \ref{Fig1}.
\setcounter{figure}{0}
  \begin{figure}[ht]
	\centering
	\includegraphics[scale=0.3]{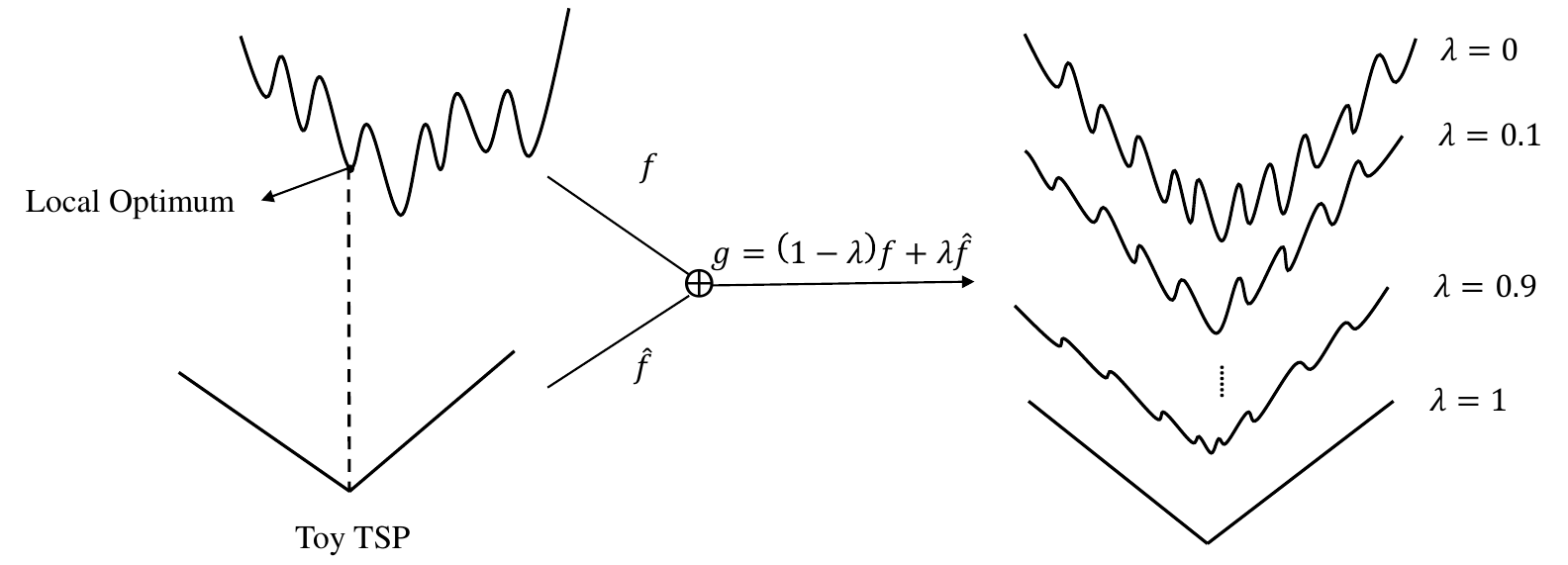}
	\caption{Effect of the proposed HC transformation for TSPs.}
	\label{Fig1}
\end{figure}

We generalize the landscape smoothing technique based on HC transformation to the UBQP. To achieve this goal, the key is how to build a unimodal smooth UBQP. For a given UBQP instance, we first construct a same size toy UBQP in which the matrix $Q$ is formed by -1 and 1 based on a known locally optimal solution of the original UBQP. We theoretically prove that this toy UBQP instance is unimodal, signifying that regardless of the initial solutions, a local search process will consistently converge towards the sole global optimum of this toy UBQP, namely, the known local optimum of the original UBQP. Then, the fitness landscape of the original UBQP can be smoothed through the application of a convex combination between the constructed toy UBQP and the original UBQP. This combination is governed by a coefficient denoted as $\lambda$, which ranges from 0 to 1. When $\lambda = 0$, the smoothed UBQP is identical to the original UBQP, preserving its original characteristics. Conversely, when $\lambda = 1$, the smoothed UBQP corresponds to the constructed unimodal toy UBQP. This process effectively defines the HC transformation, establishing a path from the original UBQP to the unimodal toy UBQP.

We subsequently present an algorithm called LSILS that integrates the landscape smoothing technique with the well-known iterated local search (ILS) \citep{lourencco2019iterated}. LSILS applies the proposed HC transformation method to smooth the landscape of the target UBQP instance and searches on the smoothed landscape. LSILS iteratively updates the smoothed UBQP landscape to reflect the most recent best found solution in relation to the original UBQP. In the experimental study, LSILS is compared against the original ILS and GH \citep{gu1994efficient} on 10 UBQP instances. GH is a  widely used method to smooth the landscape of a TSP by edge cost manipulation. Results show that HC transformation can outperform its counterparts.

Parallel metaheuristics have drawn a lot of attention from researchers with the advancement of computer hardware, especially in the context of tackling challenging optimization problems that call for sizable computational resources. To fully make use of multi-core processors to speed up metaheuristics and improve the quality of the solutions, we further implement the proposed LSILS in a parallel cooperative algorithmic framework, resulting in PC-LSILS. Processes share their best found solutions to their neighbor processes, meanwhile receive new best solutions from neighbors if they are better than the solution obtained by local search. Experimental studies on both UBQPs and TSPs demonstrate that, by utilizing parallelization, it is possible to further enhance the landscape smoothing effect of the HC transformation.

The rest of this paper is organized as follows. Section \ref{2} presents the related work. Section \ref{3} presents the construction of the toy UBQP and gives a proof that its landscape is unimodal. In Section \ref{4}, the algorithm LSILS is presented for UBQPs and then Section \ref{5} executes LSILS in parallel for both UBQPs and TSPs. The experimental results are also provided in these two sections. Section \ref{6} concludes this article and discusses future work.

\section{Related Work}\label{2}
The UBQP problem can be formalized as follows: % \citep{liefooghe2014hybrid},
\begin{equation}
	\begin{aligned}
		&\textrm{maximize} \quad f(x) = x^TQx = \sum_{i=1}^{n}\sum_{j=1}^{n}q_{ij}x_ix_j \\
		&\textrm{subject to}  \qquad \qquad x\in \{0,1\}^n
	\end{aligned}
	\label{ubqp}
\end{equation}
where $Q = [q_{ij}]$ is an $n$-dimensional matrix and $x$ is an $n$- dimensional vector with binary variables. In this article, we focus on symmetric UBQPs, meaning $q_{ij} = q_{ji}$ for all $i,j\in\{1,2,...,n\}$.

The TSP is described as follows. The objective is to identify a tour of cities that minimizes the total length of the cycle while ensuring that each city is visited only once and returns to the starting city. Let $G = (V,E)$ be a complete graph, where $V$ is a set of vertices and $E$ is a set of edges connecting the vertices.  Each edge $(i, j)$ presents a non-negative distance associated with $E$, denoted as $d_{ij}$. The TSP can be written as
\begin{align}
	&\textrm{minimize}\quad D(x)=d_{x_{(n)}x_{(1)}}+\sum_{i=1}^{n-1}d_{x_{(i)}x_{(i+1)}}\nonumber\\
	&\textrm{subject to}\quad x=\left(x_{(1)},x_{(2)},...,x_{(n)}\right)\in\mathcal{P}_n
\end{align}
where $D:\mathcal{P}_n\rightarrow\mathbb{R}$ is the objective function and $\mathcal{P}_n$ is the permutation space of $\{1, 2,...,n\}$. This paper mainly focuses on symmetric TSPs, where the distance between two cities has $d_{ij} = d_{ji}$.

Both the UBQP and the TSP can be classified into the class of NP-hard combinatorial optimization problems, which have rugged landscapes. Landscape characteristics can significantly impact the efficiency and effectiveness of search algorithms applied to UBQPs and TSPs.

\subsection{Fitness Landscape Analysis}
Nowadays, fitness landscape analysis is employed in a variety of problems to better understand how search algorithms behave. \cite{tari2018sampled} used sampled walk (SW) algorithm, a local search based on randomly sampled neighborhoods to investigate the landscape of UBQPs and found UBQP landscapes are locally rugged but globally smooth, with local optima close to each other. These findings would suggest a big-valley structure in these landscapes, which was also observed in the same year \citep{tari2018worst}. It gave an explanation that the reduced diversity of the best trajectory peaks was caused by the big valley structure of the considered UBQP landscapes, meaning local optima tend to be grouped in a single part of the landscapes. Not only UBQPs, many landscapes derived from specific combinatorial problems tend to have a big valley structure, where the local optima are clustered around a central global optimum. For example, \cite{boese1995cost} first observed a big valley structure of a TSP instance' s (att532) landscape which manifests a strong correlation between the distance of a solution to the global optimum and its cost. \cite{hains2011revisiting}  and \cite{shi2018eb} confirmed the existence of the big valley structure on some TSP instances. Based on features of fitness landscape in the previous related researches, Merz and Katayama \citeyearpar{merz2004memetic} presented a memetic algorithm and found that there was high correlation between fitness and distance to the best-known or optimum solution. Chicano and Alba \citeyearpar{chicano2013elementary} proved that the UBQP could be written as a sum of two elementary landscapes and such decomposition could be useful to create new search methods for the UBQP. \cite{tariUseLambdaEvolution2021a} deeply investigated the sampled walk search algorithm and proposed a preliminary hybridization between SW and a local search. Results highlighted the potential of combining SW with advanced intensification processes on NK landscapes but was not efficient on UBQP landscapes, so other strategies needed to be applied to exploit SW best trajectory peaks. In the work by Tari and Ochoa \citeyearpar{tari2021local}, they used a recently proposed method, the monotonic LON (MLON) model \citep{ochoa2017understanding}, to study the induced fitness landscapes to understand the performance differences among five different pivoting rules.

\subsection{Landscape Smoothing Methods}
The challenge of COPs comes from their rugged landscapes, so several studies aim to enhance the efficiency of search algorithms by smoothing the landscapes. For instance, for TSPs, Gu and Huang \citeyearpar{gu1994efficient} proposed a landscape smoothing method by edge cost manipulation. Specifically, for a given TSP, let $n$ be the number of cities, $d_{ij}$ be the distance between city $i$ and city $j, (i, j = 1, 2, . . . , n)$. With the parameter $\alpha$, the smoothed distance can be defined as
\begin{equation}
	d_{ij}(\alpha) = \left\{
	\begin{aligned}
		\overline{d}+(d_{ij}-\overline{d})^\alpha && {\text{if}\hspace{0.5em}  d_{ij}\geq\overline{d}} \\
		\overline{d}-(d_{ij}-\overline{d})^\alpha && {\text{if} \hspace{0.5em}  d_{ij}<\overline{d}} \\
	\end{aligned}	
	\right.
	\label{dij}
\end{equation}
with
\begin{equation}
	\overline{d}=\frac{1}{n(n-1)}\sum_{i,j=1, i\neq j}^n d_{ij},
\end{equation}
where $\alpha\geq 1$ is the smoothing factor and decreases step by step. $\overline{d}$ is the average distance of all cities. If $\alpha=1$, the smoothed problem is the original TSP. On the other hand, if $\alpha\to \infty, d_{ij}(\alpha)\to \overline{d}$, meaning all solutions have the same fitness, that is, the landscape is smoothed to a plane. By changing $\alpha$, a series of smoothed TSPs can be generated, and local search methods can be used for solving these problems. Based on Eq.~(\ref{dij}), \cite{schneider1997search} proposed several new smoothing functions for the TSP (linear smoothing, power-law smoothing, exponential smoothing, hyperbolic smoothing, sigmoidal smoothing and logarithmic smoothing) and combined them with simulated annealing (SA) and great deluge algorithm (GDA). \cite{coy1998see} conducted experiments to illustrate how the smoothing algorithm works, revealing that the algorithm can escape from local minima by moving uphill occasionally and by avoiding downhill movements. \cite{liang1999combining} used a local search operator to smooth the fitness landscape and prevent trapping in local optima. One year later, \cite{coy2000computational} proposed the sequential smoothing algorithm (SSA) that alternates between convex and concave smoothing function to avoid being trapped in poor local optima. \cite{dong2006stochastic} introduced the idea of probabilistic acceptance inspired by SA and developed a new algorithm that employs exponential smoothing functions for TSPs. Hasegawa and Hiramatsu \citeyearpar{hasegawa2013mutually} suggested that metropolis algorithm (MA) can be used effectively as a local search algorithm in search-space smoothing strategies and proposed a new landscape smoothing algorithm called MASSS. Recently, \cite{shi2020homotopic} proposed a new landscape smoothing method called HC transformation for the TSP. The smoothed TSP is an HC combination of the original TSP and a convex-hull TSP which is well designed at first.

To the best of our knowledge, there are no global search approaches based on landscape smoothing for the UBQP. In this paper, we take inspiration from \cite{shi2020homotopic} mentioned above and apply this HC transformation to the UBQP.

\subsection{Parallel Metaheuristics}

With the growing computational power requirements of large scale applications, parallel computing has drawn increasing interest especially in solving complex optimization problems that require large amounts of computational power. Multi-core and distributed metaheuristics are effective ways to speed up local search methods and improve the quality of the solutions in the design of more robust algorithms by using several processing units simultaneously. \cite{alba2013parallel} described a taxonomy for parallelization models of trajectory-based metaheuristics including parallel multistart, parallel moves, and move acceleration.

For ILS algorithm, \cite{araujo2007exploring} presented four slightly different strategies (independent processes, one-off cooperation, one elite solution, a pool of elite solutions) for parallelizing an extended greedy randomized adaptive search procedures (GRASPs) with iterated local search (ILS) heuristic for the mirrored TSP. Similarly, Cordeau and Maischberger \citeyearpar{cordeau2012parallel} proposed a distributed parallel tabu search (TS) with ILS framework for vehicle routing problems (VRPs) and several of its variants, while Rocki and Suda \citeyearpar{rocki2012large} proposed a cooperative parallel algorithm embedded 2-opt local search within ILS framework for the TSP. \cite{juan2014using} proposed a parallel ILS-based algorithm for permutation flow-shop problems (PFSPs). Beirigo and dos Santos \citeyearpar{beirigo2015parallel} proposed several parallelization tweaks on ILS for travel planning problems. Additionally, in parallel variable neighborhood search, \cite{yazdaniFJSP2010} proposed a parallel VNS algorithm that solves the FJSP to minimize makespan time based on the application of multiple independent searches. \cite{de2017parallel} introduced a cooperative parallel heuristic for the uncapacitated single allocation hub location problems, which combined multiple parallel ILS-RVND local search with path-relinking and cooperated through asynchronously exchanging solutions in a shared pool. \cite{herran2020parallel} developed new local search procedures integrated into the basic VNS (BVNS), and the best VNS variant was parallelized for the obnoxious p-median problem. \cite{akbay2020parallel} proposed a parallel VNS combined with quadratic programming to develop an efficient solution approach for cardinality constrained portfolio optimization. \cite{zang2022parallel} proposed a synchronous “master–slave” parallel VNS (SMPVNS) and an asynchronous cooperative parallel VNS (ACPVNS) to solve bilevel covering salesman problems (BCSPs). For tabu search (TS) and guided local search (GLS), \cite{chang2021parallel} proposed an effective parallel iterative solution-based tabu search algorithm to solve obnoxious p-median problems using a delete-add compound move strategy. \cite{uchronski2021parallel} proposed a parallel tabu search algorithm with blocks elimination criteria for weighted tardiness single-machine scheduling problems. Alalmaee and Tairan \citeyearpar{alalmaee2022hybrid} suggested a hybrid cooperative algorithm of GLS for TSPs, preserving the good components of each process by utilizing the feature of the crossover operators at the exchange information point. As for SA, \cite{sheikholeslami2021power} focused on Fujitsu’s Digital Annealer Unit (DAU) to explore the power of parallelism in the SA for the global optimum in high-dimensional optimization problems, including UBQPs, quadratic assignment problems (QAPs) and so on. More information about DA can be seen in \citep{aramon2019physics, matsubara2020digital}.

\section{HC Transformation for UBQPs}\label{3}
In this section, we propose the HC transformation for the UBQP to smooth its rugged landscape while maintain the information in elite solutions. The HC transformation changes the original UBQP by combining it with a well-designed toy UBQP. In the following, we will first construct the toy UBQP, and then give a proof that the designed toy UBQP is unimodal.

\subsection{Construction of Unimodal UBQPs}\label{AA}
The fitness landscape of the UBQP  is characterized by a big valley structure with numerous local optima. Our goal is to smooth the rugged landscape while preserving the big valley structure. Given a solution $x$, if we can construct a unimodal problem with a unique extreme point, then the landscape of this problem will manifest the effect we desire. It is worth noting that the distinguishing feature of the HC transformation, as opposed to previous landscape methods, is its capability to maintain the optimal solutions following the smoothing process.

Here, we show that for an $n$-dimensional UBQP, one can construct an $n$-dimensional unimodal UBQP based on a given solution $x$. We prove that in the unimodal UBQP, $x$ is the only $k-Bit$ locally optimal solution for $k = 1,2,\cdots,n$, meaning that regardless of the initial solutions, a bit-flip based local search process will consistently converge to $x$. The construction process is as follows.

Given a solution $x$ with size $n$, a special matrix $\hat{Q}$ can be created as
\begin{equation}
	\label{Q}
	\hat{Q}_{ij} = \left\{
			\begin{aligned}
				1,\quad
				\text{if}  \hspace{0.5em} x_ix_j=1\\
				-1,\quad \text{if}   \hspace{0.5em} x_ix_j\neq1
			\end{aligned}
		\right.,\ \ \  i,j\in\{1,2,\cdots,n\}
\end{equation}
where $x_i$ and $x_j$ denote the $i$-th and $j$-th binary variable (0-1) of the solution $x$. $\hat{Q} = [\hat{Q}_{ij}]$ is an $n \times n$ matrix of the unimodal UBQP we construct.

For example, for a 5-dimensional UBQP, given a solution $x = (0,1,0,1,1)$. We construct a UBQP that is unimodal and has $x$ as the global optimum. Initially, it can be seen that $x_1=0$, $x_2=1$, $x_3=0$, $x_4=1$, $x_5=1$, according to Eq.~(\ref{Q}), we have $\hat{Q}_{22} = \hat{Q}_{44} = \hat{Q}_{55} = \hat{Q}_{24} = \hat{Q}_{42} = \hat{Q}_{25} = \hat{Q}_{52} =\hat{Q}_{45} = \hat{Q}_{54} = 1$, the rest elements of $\hat{Q}$ are all equal to $-1$, that is
\[
\hat{Q} = \begin{pmatrix}
	-1 & -1 & -1 & -1 & -1 \\
	-1 & 1 & -1 & 1 & 1 \\
	-1 & -1 & -1 & -1 & -1 \\
	-1 & 1 & -1 & 1 & 1 \\
	-1 & 1 & -1 & 1 & 1 \\
\end{pmatrix}
\]

Based on the above design, the unimodal UBQP is constructed successfully. To prove that this designed problem is unimodal, it is necessary to show that $x$ is the unique $k-Bit$ ($\forall k \in \{1,2,\cdots,n\}$) optimal solution. The strict theoretical proof is as follows.

\subsection{Proof of the Unimodal Property}
First, we prove that in the unimodal UBQP defined by Eq.~(\ref{Q}), $x$ is a $k-Bit$ optimal solution for any $k \in \{1,2,\cdots,n\}$, then prove its uniqueness.

\noindent \textit{Definition}: A solution is said to be $k-Bit$ optimal if it is impossible to obtain another solution with a bigger objective function value by flipping any $k$ variables.

Assume $x'$ is another solution of the unimodal UBQP defined by Eq.~(\ref{Q}), which has $k$ different bits to $x$. Let $A_1 = \{n_1,n_2,\cdots,n_p\}$ denotes the set of indexes that satisfy $x_i=0$ and $x'_i=1$ and $A_2 = \{n_{p+1},\cdots,n_k\}$ denotes the index set that satisfy $x_i=1$ and $x'_i=0$, i.e. $A_1 = \{i \vert x_i=0, x'_i = 1\}$ and $A_2 = \{i \vert x_i=1, x'_i = 0\}$. Then we can use $x$ to rewrite $x'$ as
\begin{equation}\label{x_diff}
	x'=x+(\varepsilon_{n_1}+\cdots+\varepsilon_{n_p})-(\varepsilon_{n_{p+1}}+\cdots+\varepsilon_{n_k}), 0\leq p\leq k (k\geq 1)
\end{equation}
where $\varepsilon_i$ denotes an $n$-dimensional binary vector with the $i$-th element being $1$ and the other $n-1$ elements being 0.

Then we have
\begin{equation}
	\begin{aligned}
		&x'^T\hat{Q}x'\\
		= &\left(x + \sum_{i=n_1}^{n_p} \varepsilon_i - \sum_{i=n_{p\text{+1}}}^{n_k} \varepsilon_i \right)^T \hat{Q} \left( x + \sum_{i=n_1}^{n_p} \varepsilon_i - \sum_{i=n_{p\text{+1}}}^{n_k} \varepsilon_i \right) \\
		= &  x^T \hat{Q}x + \sum_{i,j=n_1}^{n_p} \varepsilon_i^T \hat{Q} \varepsilon_j +  \sum_{i,j=n_{p+1}}^{n_k} \varepsilon_i^T \hat{Q} \varepsilon_j - \\
		&2 \sum_{i=n_1}^{n_p} \sum_{j=n_{p+1}}^{n_k} \varepsilon_i^T \hat{Q} \varepsilon_j + 2 \sum_{i=n_1}^{n_p} x^T \hat{Q} \varepsilon_i - 2 \sum_{i=n_{p+1}}^{n_k} x^T \hat{Q} \varepsilon_i.
	\end{aligned}
\end{equation}
By the properties of the matrix $\hat{Q}$ and the vector $x$, we have
\begin{subequations}
	\label{eq9}
	\begin{align}
		x'^T&\hat{Q}x'-x^T\hat{Q}x ~ =2\sum_{i=n_1}^{n_p} x^T \hat{Q}\varepsilon_i  \\
		&-2\sum_{i=n_{p+1}}^{n_k} x^T \hat{Q}\varepsilon_i  \\
		&+\sum_{i,j=n_1}^{n_p}  \hat{Q}_{ij} +  \sum_{i,j=n_{p+1}}^{n_k} \hat{Q}_{ij} -2\sum_{i=n_1}^{n_p} \sum_{j=n_{p\text{+1}}}^{n_k} \hat{Q}_{ij}.
	\end{align}
\end{subequations}
For convenience, we divide Eq.~(\ref{eq9}) into three parts (8a), (8b), and (8c) and discuss them one by one.

We can summarize the relationship between $\hat{Q}_{ij}$ and sets $A_1, A_2$ as \tablename~\ref{A1A2}.
\begin{table}[htbp]
	\begin{center}
		\caption{The value of $\hat{Q}_{ij}$ when indexes belong to different sets.}
		\label{A1A2}
		\begin{tabular}{|c|c|c|}
			\hline
			\diagbox{i}{j}		&  $j\in A_1$                      & $j\in A_2$ \\ \hline
			$i\in A_1$ & \multicolumn{1}{c|}{$\hat{Q}_{ij}=-1$} & $\hat{Q}_{ij}=-1$ \\ \hline
			$i\in A_2$ & \multicolumn{1}{c|}{$\hat{Q}_{ij}=-1$} & $\hat{Q}_{ij}=1$  \\ \hline
		\end{tabular}
	\end{center}
\end{table}

For (8a) of Eq.~(\ref{eq9}) where $i\in A_1$, from \tablename~\ref{A1A2}, we have
\begin{equation}
\hat{Q}_{ij} = -1.
\end{equation}
$\hat{Q}\varepsilon_i$ is the $i$th column of $\hat{Q}$, $x$ is a 0-1 vector, and by assumption, at least $(k-p)$ elements equal to $1$, then for $\forall i \in A_1$, we get
$x^T \hat{Q}\varepsilon_i \leq -(k-p)$, thus
\begin{equation}
	\text{(8a)} \leq -2p(k-p).
\end{equation}

For (8b) of Eq.~(\ref{eq9}), there is $\forall i \in A_2$,
\begin{equation}
	x^T \hat{Q}\varepsilon_i = \sum_{j \in A_2} \hat{Q}_{ij} x_j + \sum_{j \notin A_2} \hat{Q}_{ij} x_j.
\end{equation}
When $i\in A_2$,
\begin{enumerate}[i)]
	\item For $j \in A_2 $, $x_j = 1$, we have $\hat{Q}_{ij}=1$, then $\sum_{j \in A_2} \hat{Q}_{ij} x_j = k-p$.
	\item For $j \notin A_2$, $x_j$ may be $0$ or $1$, for $x_j = 0$, $\hat{Q}_{ij} x_j = 0$ , and for $x_j = 1$, $\hat{Q}_{ij} x_j = 1$ .
\end{enumerate}
Thus
\begin{equation}
	x^T \hat{Q}\varepsilon_i = \sum_{j \in A_2} \hat{Q}_{ij} x_j + \sum_{j \notin A_2} \hat{Q}_{ij} x_j \geq k-p,
\end{equation}
thereby
\begin{equation}\label{7b_res}
	\text{(8b)}\leq -2(k-p)^2.
\end{equation}

For (8c) of Eq.~(\ref{eq9}),
\begin{equation}
	\text{(8c)} = -p^2 + (k - p)^2 + 2p(k - p) = k^2 - 2p^2.
\end{equation}
And we can obtain:
\begin{equation}
	\text{(8a)} + \text{(8b)} + \text{(8c)} \leq -{(k - p)}^2 - p^2 < 0,
\end{equation}
Therefore, the inequality
\begin{equation}\label{global}
	x'^T \hat{Q}x' < x^T  \hat{Q}x
\end{equation}
is proved. Combining the above analyses, we can get the inequality (\ref{global}) holds for any $x'$, meaning $x$ is a $k-Bit$ optimal solution for any $k\in\{1,2,\cdots,n\}$. Hence $x$ is a global optimal solution of the UBQP problem defined by $\hat{Q}$.

In the following, we will prove that $x$ is the unique $k-Bit$ optimal solution. Suppose that there exists another $k-Bit$ optimal solution $x'$, which has $q$ different bits to $x$. Since $x$ is $k-Bit$ optimal, it must be $k'-Bit$ optimal for $k'<k$. Hence, to prove the uniqueness of $x$, it is only necessary to prove that $x'$ is not $1-Bit$ optimal. Let
$x'=x+(\varepsilon_{n_1}+\cdots+\varepsilon_{n_w})-(\varepsilon_{n_{w+1}}+\cdots+\varepsilon_{n_q}), $
where $A'_1= \{i \vert x_i=0, x'_i = 1\}=\{n_1,n_2,...n_w\}$, $A'_2= \{i \vert x_i=1, x'_i = 0\}=\{n_{w+1},n_{p+2},...n_q\}$.
For $x''\in N_{1-Bit}(x')$, where $N_{1-Bit}(x')$ denotes the $1-Bit$ neighborhood of $x'$ that contains solutions that have only one bit opposite to $x'$, consider the following two scenarios:

\begin{enumerate}[i)  ]
	\item First, when $w\neq 0$ (i.e. $A'_1\neq \emptyset$), take an index $n_m \in A'_1$ and take $x'' = x'-\varepsilon_{n_m} \in N_{1-Bit}(x')$, we have
	\begin{equation}	
		\begin{aligned}
			x''^T \hat{Q}x''=&\left(x' - \varepsilon_{n_m} \right)^T  \hat{Q} \left( x' - \varepsilon_{n_m} \right) \\
			= & x'^T \hat{Q}x' -\varepsilon_{n_m}^T \hat{Q}x'\\
			&-x'^T \hat{Q}\varepsilon_{n_m}+\varepsilon_{n_m}^T \hat{Q}\varepsilon_{n_m},
		\end{aligned}
	\end{equation}
	where $ \hat{Q}$ is the symmetric matrix we construct, so $\varepsilon_{n_m}^T \hat{Q}x'=x'^T \hat{Q}\varepsilon_{n_m}$ and $\varepsilon_{n_m}^T \hat{Q}\varepsilon_{n_m}= \hat{Q}_{n_m,n_m}$. Since $n_m \in A'_1$, there is $x_{n_m}=0$. Then it is clear that there is $\varepsilon_{n_m}^T \hat{Q}\varepsilon_{n_m}= \hat{Q}_{n_m,n_m}=-1$. On the other hand, from $x_{n_m}=0$, we can get for any $j\in\{1,2,... .n\}$, $\hat{Q}_{n_m,j} = -1$. Then we have $\hat{Q}_{n_m,j'}x'_{j'}=-1$, for $\{j'|x'_{j'}=1\}$ and $\hat{Q}_{n_m,j''}x'_{j''}=0$, for $\{j''|x'_{j''}=0\}$. Since $x'\neq \textbf{0}$, there exists at least one variable in $x'$ equals to $1$, then
	\begin{equation}
		\varepsilon_{n_m}^T\hat{Q}x'=\sum_j\hat{Q}_{n_m,j}x'_{j}\leq \hat{Q}_{n_m,j'}x'_{j'}=-1.
	\end{equation}
	Hence we have
	\begin{equation}
		x''^T\hat{Q}x''\geq x'^T \hat{Q}x'+1>x'^T \hat{Q}x',
	\end{equation}
	which means $x'$ is not $1-Bit$ optimal under the scenario $w\neq 0$.
\item Second, when $w=0$ (i.e. $A'_1=\emptyset$), take an index $n_m\in A'_2$ and take $x'' = x'+\varepsilon_{n_m} \in N_{1-Bit}(x')$, we have
\begin{equation}
	\begin{aligned}
		x''^T\hat{Q}x''=&\left(x' + \varepsilon_{n_m} \right)^T \hat{Q} \left( x' + \varepsilon_{n_m} \right)\\
		=& x'^T\hat{Q}x' +\varepsilon_{n_m}^T\hat{Q}x'\\
		&+x'^T\hat{Q}\varepsilon_{n_m}+\varepsilon_{n_m}^T\hat{Q}\varepsilon_{n_m}.
	\end{aligned}
\end{equation}
Similarly there is $\varepsilon_{n_m}^T\hat{Q}x'=x'^T\hat{Q}\varepsilon_{n_m}$ and $\varepsilon_{n_m}^T\hat{Q}\varepsilon_{n_m}=\hat{Q}_{n_m,n_m}$, and by the assumption that there is $x_{n_m}=1$, it is clear that $\varepsilon_{n_m}^T\hat{Q}\varepsilon_{n_m}=\hat{Q}_{n_m,n_m}=1$.

\begin{equation}
	\varepsilon_{n_m}^T\hat{Q}x'=\sum_j\hat{Q}_{n_m,j}x'_j=\sum_{j \in A'_2}\hat{Q}_{n_m,j}x'_j+\sum_{j\notin A'_2}\hat{Q}_{n_m,j}x'_j,
\end{equation}

If $j\in A'_2$, $x_j = 1, x'_j = 0, \sum_{j\in A'_2}\hat{Q}_{n_m,j}x'_j=0$;

If $j\notin A'_2$, $x'_j$ may be 0 or 1, if $x'_j=1$, then $\hat{Q}_{n_m,j}x'_j=1$, and if $x'_j=0$, then $\hat{Q}_{n_m,j}x'_j=0$, so $\sum_{j\notin A'_2}\hat{Q}_{n_m,j}x'_j\geq 0$, thus $\varepsilon_{n_m}^T\hat{Q}x'\geq 0$. Hence we have
\begin{equation}
 x''^T\hat{Q}x''\geq x'^T \hat{Q}x'+1 >x'^T \hat{Q}x',
\end{equation}
which means $x'$ is not $1-Bit$ optimal under the scenario $w= 0$.
\end{enumerate}

In summary, for any $x'\neq x$, there exists $x''\in N_{1-Bit}(x')$ such that $f(x'')>f(x')$, i.e., $x'$ is not $1-Bit$ optimal. This means $x'$ cannot be a $k-Bit$ optimal solution. Thus we prove that $x$ is the unique $k-Bit$ optimal solution to the well-designed UBQP problem for $\forall k\geq1$ $\Box$.

\subsection{HC Transformation of UBQPs}
Based on the above construction steps, we obtain a well-designed matrix, then for a UBQP, when an elite solution, denoted as $x_e$, is obtained, we can create a unimodal UBQP based on $x_e$. By combining the unimodal UBQP with the original UBQP, the solution space of the original UBQP can be smoothed and at the same time, the information contained in $x_e$ can be preserved. The HC transformation of UBQP is defined by a convex combination of the unimodal UBQP and the original UBQP. The specific combining process is as follows. In the smoothed UBQP, the matrix $Q'=[Q'_{ij}]$ is set to be
\begin{equation}
	Q'_{ij}(\lambda) = (1-\lambda)Q_{ij} + \lambda \hat{Q}_{ij},
\end{equation}
where $Q = [Q_{ij}]$ is the matrix of the original UBQP and $\hat{Q} = [\hat{Q}_{ij}]$ is the matrix of the unimodal UBQP. The objective function after HC transformation ($g$) is expressed as
\begin{equation}
	g(x|\lambda) = (1-\lambda)f_o(x) + \lambda \hat{f}(x),
\end{equation}
where $f_o$ is the objective function of the original UBQP and $\hat{f}$ is the objective function of the unimodal UBQP constructed based on a local optimal solution. From the definition of $g$, it can be seen that when $\lambda=0$, the smoothed UBQP $g$ degenerates to the original UBQP $f_o$, when $\lambda=1$, it is smoothed to the unimodal UBQP $\hat{f}$. Actually, $g$ is a homotopic transformation from $f_o$ to $\hat{f}$, with this kind of transformation, the elite solution can be preserved and the landscape becomes smoother.

\section{Algorithmic Framework}\label{4}
In this section, we give a framework of LSILS for UBQPs based on the definition of HC transformation proposed in section \ref{3}. Then a parallel cooperative algorithm framework of LSILS is introduced to utilize the capabilities of multiple search processes operating in parallel.

\subsection{LSILS for UBQPs}
The entire procedure of the LSILS process is summarized in Algorithm \ref{LSILS}, it iteratively executes a local search procedure and a perturbation procedure.

\begin{algorithm}
		\caption{Landscape Smoothing Iterated Local Search for UBQPs}
		\label{LSILS}
		$x_{ini}\leftarrow$ Randomly generated solution\;
		$x_{(0)}\leftarrow$ ILS$(x_{ini}|f_o)$\;
		$x_{f_o,best}\leftarrow x_{(0)}$\;
		$j\leftarrow 0$\;
		\While{stopping criterion is not met}
		{
			Construct the unimodal UBQP $\hat{f}$ based on $x_{f_o,best}$\;
			$g\leftarrow(1-\lambda)f_o+\lambda \hat{f}$\;
			$x'_{(j)}\leftarrow$Perturbation($x_{(j)}$)\;
			$\{x_{(j+1)},x_{f_o,best}\}\leftarrow$LS($x'_{(j)},x_{f_o,best}|g$)\;
			$j\leftarrow j+1$\;
			$\lambda$ = Update($\lambda$)\;
		}
		\Return $x_{f_o,best}$
\end{algorithm}

The initial local optimum $x_{(0)}$ is obtained by ILS from a random generated solution $x_{ini}$ (line 2) and $x_{f_o,best}$ is the current best solution with respect to (w.r.t.) $f_o$ which is first assigned to $x_{(0)}$ (line 3). Then LSILS repeats the following steps until the stop condition is met (line 5-11). At each iteration, first, based on $x_{f_o,best}$, the unimodal UBQP $\hat{f}$ is constructed (line 6). The objective function of the smoothed UBQP $g$ is obtained based on $f_o$ and $\hat{f}$ with the smoothing factor $\lambda$ (line 7). The local search is then applied on the smoothed UBQP from a perturbed solution $x'_{(j)}$ (line 8) of the current solution $x_{(j)}$. LS($x'_{(j)}$, $x_{f_o,best}|g$) means executing a local search from $x'_{(j)}$ on the smoothed UBQP $g$, meanwhile keeping updating $x_{f_o,best}$ by tracking the value change of $f_o$. Its return value $x_{(j+1)}$ is the local optimum w.r.t. the smoothed UBQP $g$ (line 9). $x_{(j+1)}$ will be used as the input of the next LS after being perturbed (line 10). Finally the coefficient $\lambda$ is updated (line 11). The update method of $\lambda$ is defined by the users.

\subsection{Parallel Cooperative LSILS}
To fully make use of multi-core processes, in this paper, we propose a parallel LSILS algorithm and integrate a cooperative mechanism. This parallel cooperative LSILS algorithm incorporates inter-process communication among processes, they will share their best found solutions while conduct local search at the same time. The entire procedure of a PC-LSILS process is shown in
Algorithm \ref{P-LSILS}. Lines 6-9 correspond the communication phase, lines 10-15 correspond the LSILS phase.

Starting from a randomly generated solution, each PC-LSILS process performs an ILS procedure from a different initial solution. If the stop condition is not satisfied, the process continues running the ILS search until it does. After the local search, the process moves on to the perturbation phase to perturb the current optimal solution w.r.t. the smoothed UBQP $x_{(j)}$ and uses the perturbed solution $x'_{(j)}$ as the input for the next round of ILS.

\begin{algorithm}
	\caption{Parallel Cooperative Landscape Smoothing Iterated Local Search}
	\label{P-LSILS}
	$x_{ini}\leftarrow$ Randomly generated solution\;
	$x_{(0)}\leftarrow$ ILS$(x_{ini}|f_o)$\;
	$x_{f_o,best}\leftarrow x_{(0)}$\;
	$j\leftarrow 0$\;
	\While{stopping criterion is not met}
	{
		\If {$x_{f_o,best}$ has updated}{ SendToNeighbors($x_{f_o,best}$)\;
			$S_r\leftarrow$ TrytoReceive()\;
			$x_e\leftarrow$SelectBestSolution($S_r\cup {x_{f_o,best}}$)\;
		}
		Construct  the unimodal UBQP $\hat{f}$ based on $x_e$\;
		$g\leftarrow(1-\lambda)f_o+\lambda \hat{f}$\;
		$x'_{(j)}\leftarrow$Perturbation($x_{(j)}$)\;
		$\{x_{(j+1)},x_{f_o,best}\}\leftarrow$LS($x'_{(j)},x_e|g$)\;
		$j\leftarrow j+1$\;
		$\lambda$ = Update($\lambda$)\;
	}
	\Return $x_{f_o,best}$
\end{algorithm}

\begin{figure}[htbp]
 	\centering
 	\includegraphics[scale=0.2]{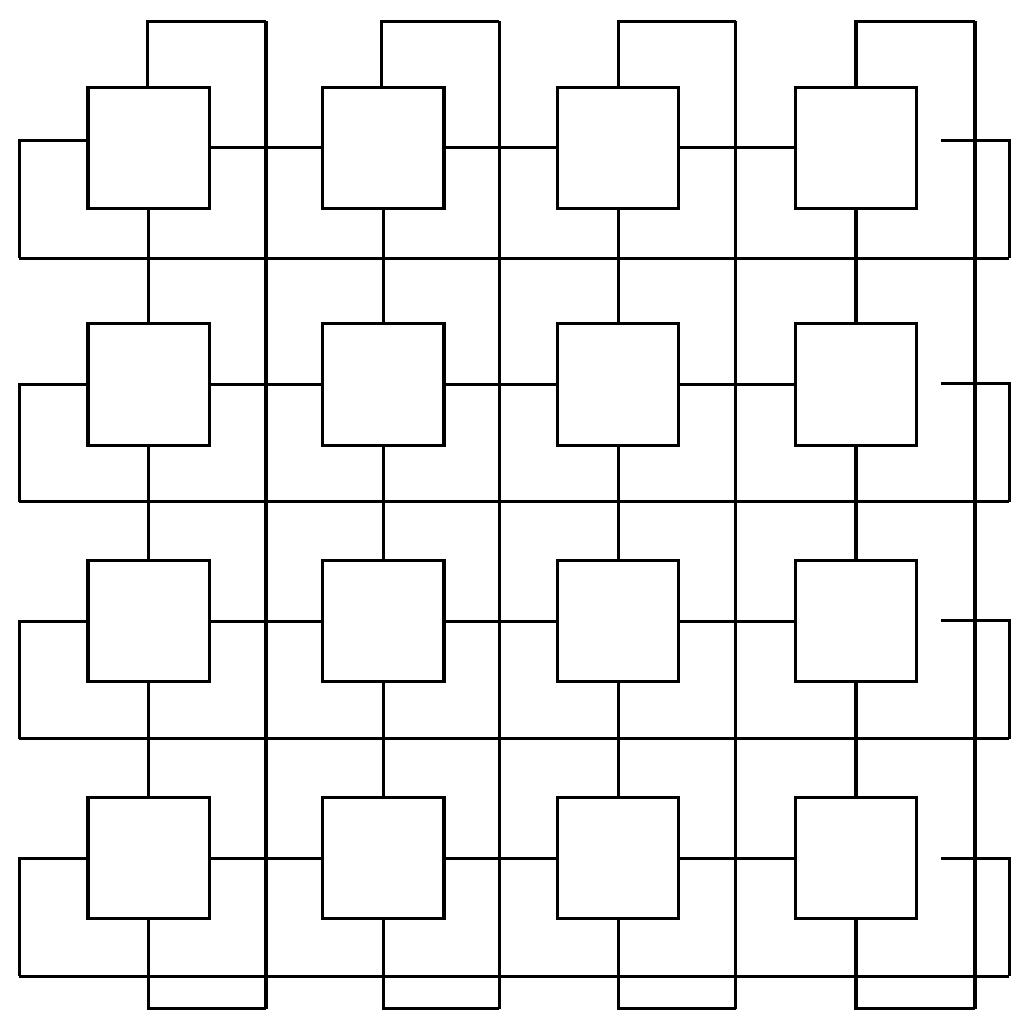}
 	\caption{A $4\times 4$ torus topology of PC-LSILS.}
 	\label{Fig2}
 \end{figure}

Since too much communication may result in having less time for local search, to reduce the communication load, the communication topology of PC-LSILS follows a torus topology, shown as \figurename~\ref{Fig2}. Each process has four neighbors and only communicates with its four neighbors. Specifically, when a new optimal solution is discovered, it is transmitted to the four neighbors (line 7). Simultaneously, processes will check if there are any new solutions from neighbors (line 9). If new solutions exist, then receive these new solutions, let $S_r$ denotes the set of solutions from neighbors. In addition, each process will also conduct its own search and obtain the historical best solution $x_{f_o,best}$ found by itself. Note that although a process may receive better solutions from its neighbors, it always sends the best solution found by itself to its neighbors. Every process maintains an elite solution $x_e$ which is the best solution in the set $\{S_r \cup x_{f_o,best}\}$. This elite solution $x_e$ will be used for the construction of the unimodal toy UBQP and guides the local search.

\section{Experimental Studies}\label{5}
\rv{In order to investigate the performance of the landscape smoothing effect of HC transformation, we first conduct landscape analysis experiments and then carry out experiment with the sequential version of LSILS algorithm.} Ten UBQP instances with a size of $n = 2500$ are chosen from ORLIB \citep{beasley1996obtaining}. These ten instances are named as {bqp2500.1, ..., bqp2500.10}, and all possess a density of 0.1, which denotes the ratio of non-zero values in the matrix $Q$. All of these ten instances have known best solutions. \rv{We conduct a landscape analysis on 5 of these 10 instances.} Our next experiment aims to parallelize the LSILS algorithm and we test the landscape smoothing effect of parallel local search with HC transformation on both UBQPs and TSPs. Ten TSP instances are selected from the TSPLIB \citep{reinelt1991tsplib} in the second experiment with the city size $n$ ranging from 400 to 1817, including rd400, u574, p654, d657, u724, rat783, pcb1173, rl1304, vm1748, and u1817.

Experiments in this article are executed on the Tianhe-2 supercomputer. Tianhe-2 is one of the world' s top-ranked supercomputers. It is equipped with 17,920 computer nodes, each comprising two Intel Xeon E5-2692 12C (2.200 GHz) processors. So each node has 24 cores and the communication between different processes is achieved using MPI (Message Passing Interface). MPI is a widely used communication protocol and library for parallel computing, allowing processes in a distributed computing environment to exchange information and coordinate their actions efficiently.

\rv{\subsection{Effects of the HC Transformation}\label{Effects}
To illustrate the smoothing effects of the HC transformation on UBQPs, in the following, we conduct landscape analysis experiments to study the landscape after the HC transformation. We execute ILS on the transformed UBQP and use the two metrics to characterize the landscape, which are \textit{local optimum density} and \textit{escaping rate}. Local optimum density is defined as the number of local optima encountered by an ILS process per move. Here a move means the ILS agent moves to a new solution from its origanl solution. Then the local optimum density can be calculated by 
\begin{equation}
    \textit{local optimum density} = \frac{N_{LO}}{Moves},
\end{equation}
where $N_{LO}$ denotes the number of local optima encountered during search and $Moves$ is the total number of moves made by the ILS algorithm. Escaping rate is defined as the success rate that a new local optimum is reached by ILS starting from a new solution obtained by perturbing the current local optimum. Then the escaping rate can be calculated by 
\begin{equation}
    \textit{escaping rate} = \frac{N_{succ}}{N_{pert}},
\end{equation}
where $N_{succ}$ is the number of successful attempts where a new local optimum is reached and $N_{pert}$ is the total number of perturbations (attempts to escape the current local optimum).
}

\rv{Local optimum density and escaping rate provide valuable insights into the ruggedness of the landscape. A higher local optimum density typically indicates a more rugged landscape, making the optimization problem more challenging to solve using local search methods. Additionally, a higher escaping rate suggests a more rugged landscape, as the algorithm frequently encounters and must escape numerous local optima. This combination of metrics can measure the difficulties posed by highly rugged landscapes in optimization tasks.}

\rv{In the following, we investigate the smoothing effects of the proposed HC transformation by using either global optimum or local optimum. First, since all test UBQP instances have known global optima, we use the global optima to construct the toy unimodal UBQP and conduct HC transformation with different $\lambda$ values. We execute $100\,000$ moves on the smoothed UBQP from a random initial solution where one move means the local search algorithm moves from the current solution to a new solution. Each sampling experiment is conducted 20 times and then averaged.}

\begin{figure*}[!t]
    \centering
	\subfloat[]{\includegraphics[scale=0.55]{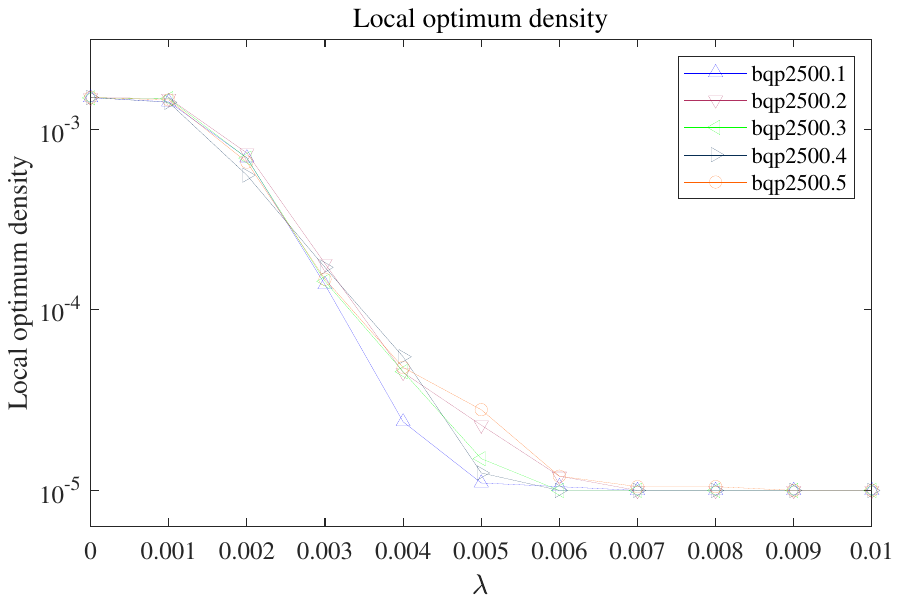}
		\label{LO_density_global}}
    \hspace{1em}
	\subfloat[]{\includegraphics[scale=0.55]{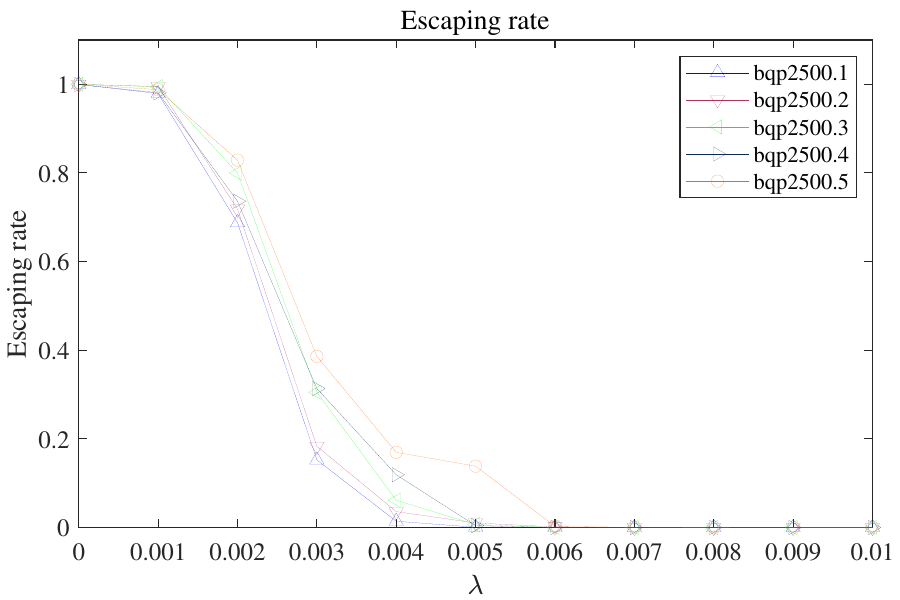}
		\label{Escaping_rate_global}}
	\caption{Landscape analysis results of the smoothed UBQPs with different $\lambda$ values based on the global optimum of the original UBQP}
    \label{LA_global}
\end{figure*}

\rv{\figurename~\ref{LO_density_global} shows how the local optimum density of the smoothed UBQP landscape changes against $\lambda$. It reflects that the local optimum density decreases as $\lambda$ increases, indicating that the HC transformation can effectively smooth the UBQP landscape. Additionally, from \figurename~\ref{Escaping_rate_global}, we can see that in all instances, the escaping rate decrease as $\lambda$ increases. Similar to the local optimum density, 
the smoothing effect is related to the parameter $\lambda$. From \figurename~\ref{LA_global}, we can conclude that the HC transformation based on global optima can indeed smooth the landscape of UBQP.}

\rv{Notice, however, that it is not always possible to know the global optimum of a problem. Therefore, we investigate the smoothing effects of the HC transformation by using a local optimum to construct the toy UBQP.}%For each instance, we first run a tabu search from a random solution to get a local optimum. In tabu search algorithm, we set $TabuT enure(i) = c + rand(10)$, where $c$ is a constant equals $n/100$ and $rand(10)$ denotes a randomly generated number from 1 to 10. Based on this obtained local optimum, then we can construct the toy UBQP with different $\lambda$ values.

\begin{figure*}[!t]
    \centering
	\subfloat[]{\includegraphics[scale=0.55]{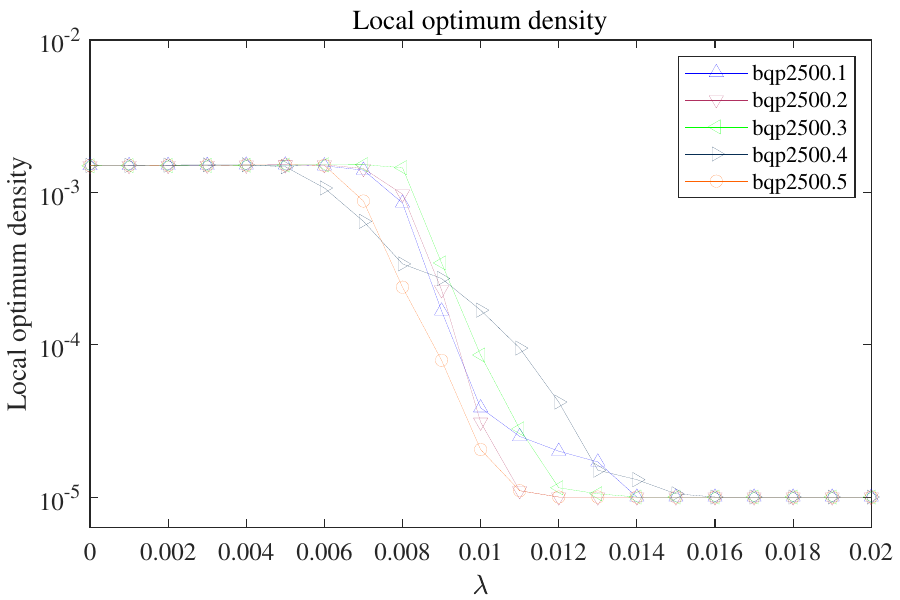}
		\label{LO_density_local}}
	\subfloat[]{\includegraphics[scale=0.55]{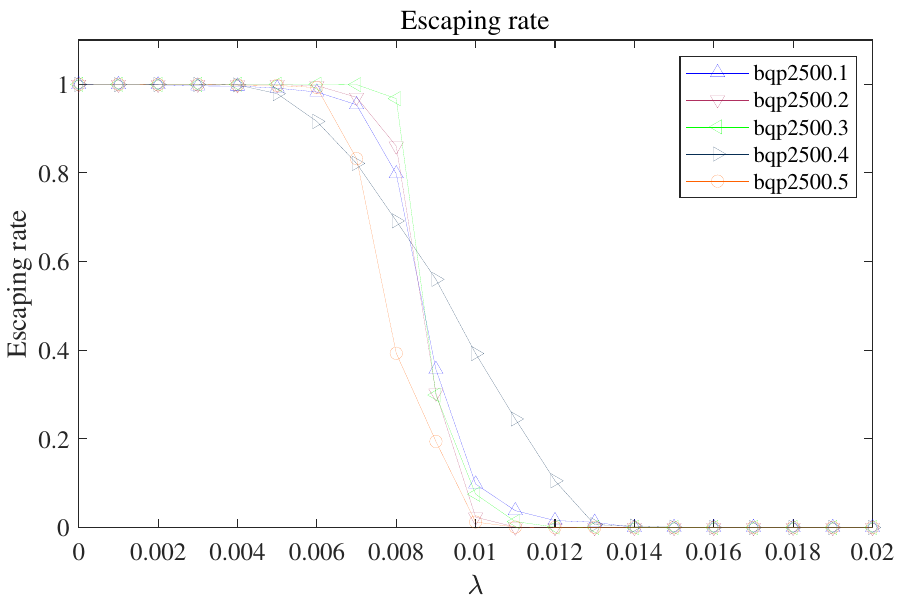}
		\label{Escaping_rate_local}}
	\caption{Landscape analysis results of the smoothed UBQPs with different $\lambda$ values based on a local optimum of the original UBQP}
    \label{LA_local}
\end{figure*}

\rv{The local optimum density values of different smoothed UBQPs are shown in \figurename~\ref{LO_density_local}. We can see that in all instances, the local optimum density is negatively related to $\lambda$. \figurename~\ref{Escaping_rate_local} shows the escaping rate versus $\lambda$ when the toy UBQP is constructed based on the local optimum. From \figurename~\ref{Escaping_rate_local} we can see that in all instances, the escaping rate is negatively related to $\lambda$. From \figurename~\ref{LA_local}, we can conclude that HC transformation based on a local optimum can also smooth the landscape of UBQP as the HC transformation based on global optima.}

\subsection{Performance of LSILS for UBQPs}\label{sectionLSILS}
\rv{Experiments in this section} aim to evaluate the smoothing effect of HC transformation on UBQPs with sequential LSILS. The developed LSILS is compared against ILS which lacks HC transformation and GH, a landscape smoothing algorithm proposed by Gu and Huang \citeyearpar{gu1994efficient}. To apply this GH method to UBQPs, based on the idea of GH, we propose a definition of the smoothed  $\widetilde{Q}$ matrix as
\begin{equation}
	\widetilde{Q}_{ij} = \left(\frac{Q_{ij}}{\left|Q_{ij}\right|_{max} + 1}\right)^\alpha,
\end{equation}
where $\alpha\geq 1$ is the smoothing factor, $\left|Q_{ij}\right|_{max}$ is the max absolute value of $Q$ matrix. When $\alpha = 1$, it is just the original UBQP with a scaling parameter $1/ \left(\left|Q_{ij}\right|_{max} + 1\right)$. From the definition of $\widetilde{Q}$, we have $\lim_{\alpha\rightarrow +\infty}\widetilde{Q}_{ij} = 0$. That is, with an increasing $\alpha$, all elements of matrix $\widetilde{Q}$ approach to a fixed value 0. This implies that all solutions will have the same fitness value, that is, the UBQP landscape is smoothed to a plane. The parameter settings are the same as \citep{gu1994efficient}, for GH, the smoothing factor $\alpha = 6 \rightarrow 5 \rightarrow 4 \rightarrow 3 \rightarrow 2 \rightarrow 1$ in the first six local search rounds. Then the algorithm will execute ILS on the original UBQP landscape ($\alpha = 1$) until the end of run time.

We have known from experiments on UBQPs that the \emph{first improvement} local search strategy and the \emph{best improvement} local search strategy have different behaviors on UBQP landscapes \citep{tari2018worst}. Based on our preliminary results, we find that in most cases, \emph{best improvement} performs better than  \emph{first improvement}, so we will apply the \emph{best improvement} strategy to UBQPs in the following experiments.

Since matrices $Q$ of UBQP instances we test are $2500 \times 2500$ with total $6.25\times 10^6$ numbers while only about $6\times10^5$ are non-zero, and they nearly range from -100 to 100, but we design the matrix $\hat{Q}$ for unimodal UBQP with numbers -1 and 1. To minimize the gap between the two matrices in order to smooth the landscape of the original UBQP without causing great damage to it, we first multiply 5 before the unimodal UBQP, getting the final smoothed $Q'$ matrix, that is
\begin{equation}
	Q'_{ij}(\lambda) = (1-\lambda)Q_{ij} + 5 * \lambda \hat{Q}_{ij}.
\end{equation}

The parameter settings are set as follows, after many preliminary attempts, we set max value $\lambda$ as 0.004, probably because we have converted a large number of zeros in the original matrix to non-zero, resulting in a relatively large change to the landscape of the UBQP. The time point to increase $\lambda$ is based on the CPU time, changing from 0 to 0.001 at 200 s, from 0.001 to 0.002 at 400 s, and so on. The stopping criterion for each algorithm is set at 1000 seconds of CPU runtime, and every 10 seconds as a time point for information logging. Each algorithm is executed for 20 independent runs on each benchmark instance. In the local search implementation, \emph{1-bit-flip} local search is used. The perturbation strength is set to be $n/4$ where $n$ is the size of instances, that means we randomly select $n/4$ variables to flip every shake.

To measure the quality of a solution, we use excess of the obtained solution, which is defined as
\begin{equation}
	\text{excess} = \frac{f(x_{\textrm{LO}})-f(x_{\text{opt}})}{f(x_{\text{opt}})},
\end{equation}
where $f(x_{\text{opt}})$ denotes the global optimum already known and $f(x_{\textrm{LO}})$ is the solution we get.

\begin{figure*}[!t]
	\subfloat[]{\includegraphics[width=1.8in]{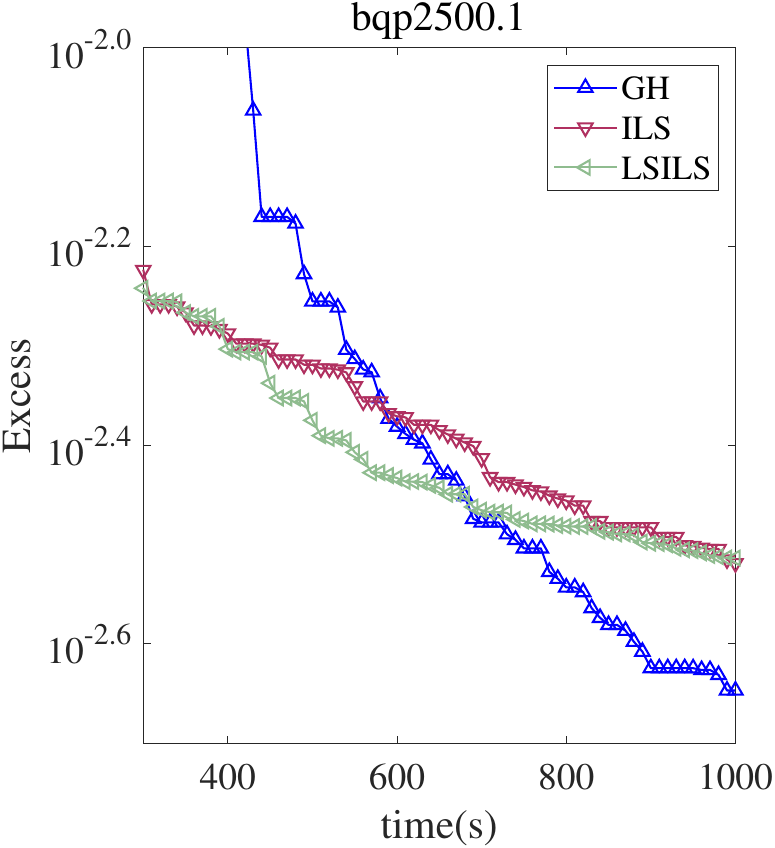}
		\label{se.bqp2500.1}}
	\subfloat[]{\includegraphics[width=1.8in]{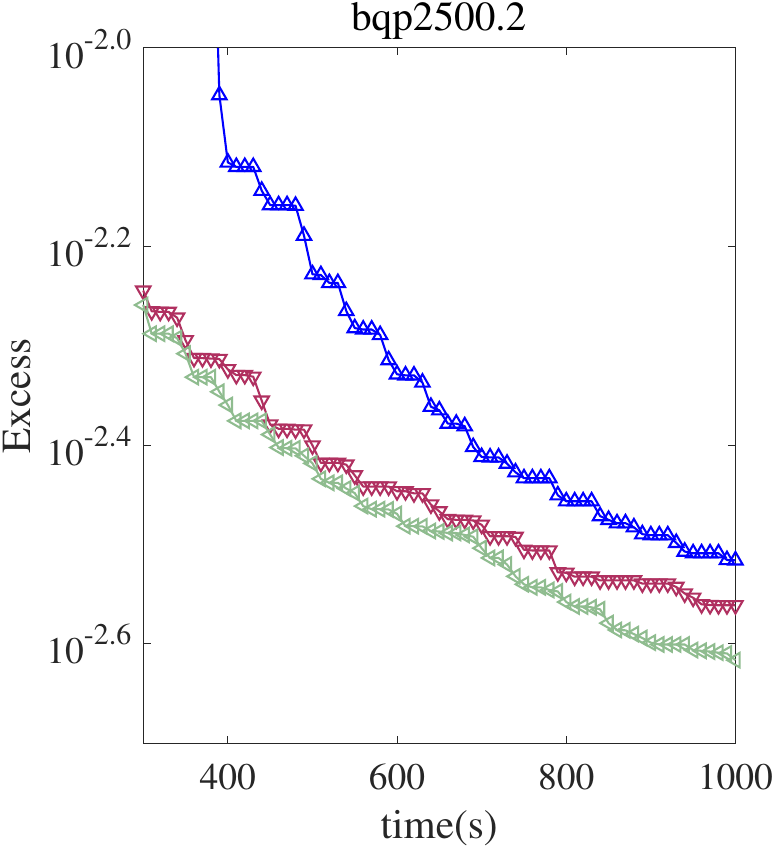}
		\label{se.bqp2500.2}}
	\subfloat[]{\includegraphics[width=1.8in]{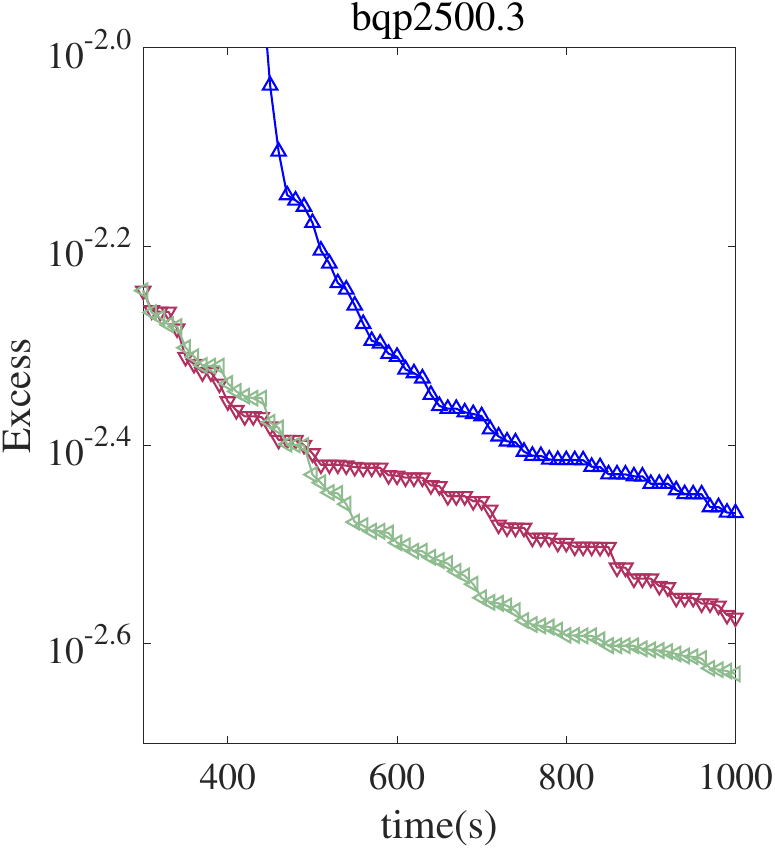}
		\label{se.bqp2500.3}}
	\subfloat[]{\includegraphics[width=1.8in]{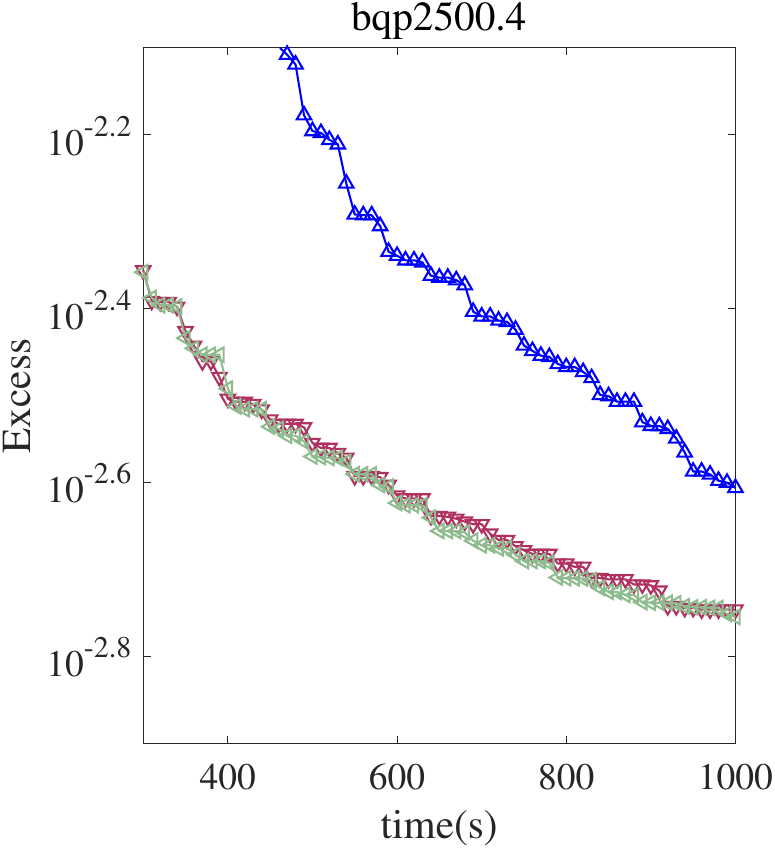}
		\label{se.bqp2500.4}}
	
	\subfloat[]{\includegraphics[width=1.8in]{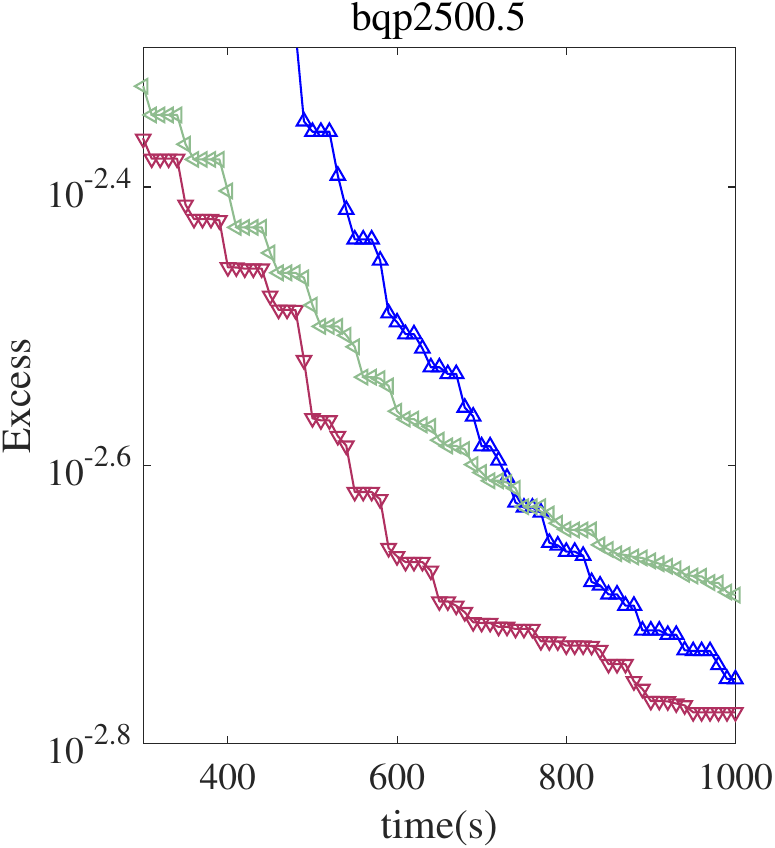}
		\label{se.bqp2500.5}}
	\subfloat[]{\includegraphics[width=1.8in]{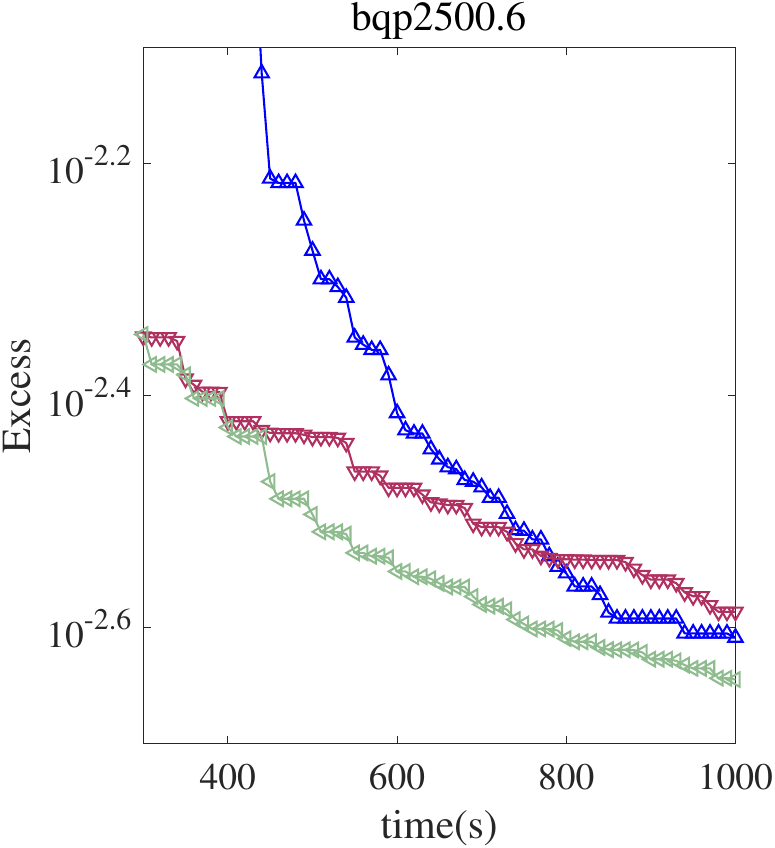}
		\label{se.bqp2500.6}}
	\subfloat[]{\includegraphics[width=1.8in]{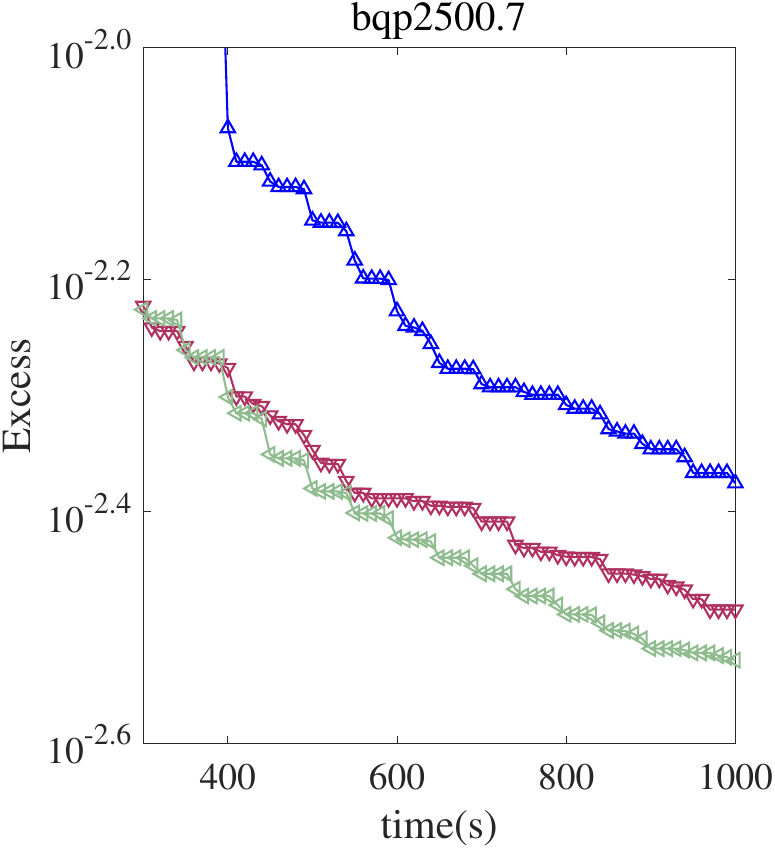}
		\label{se.bqp2500.7}}
	\subfloat[]{\includegraphics[width=1.8in]{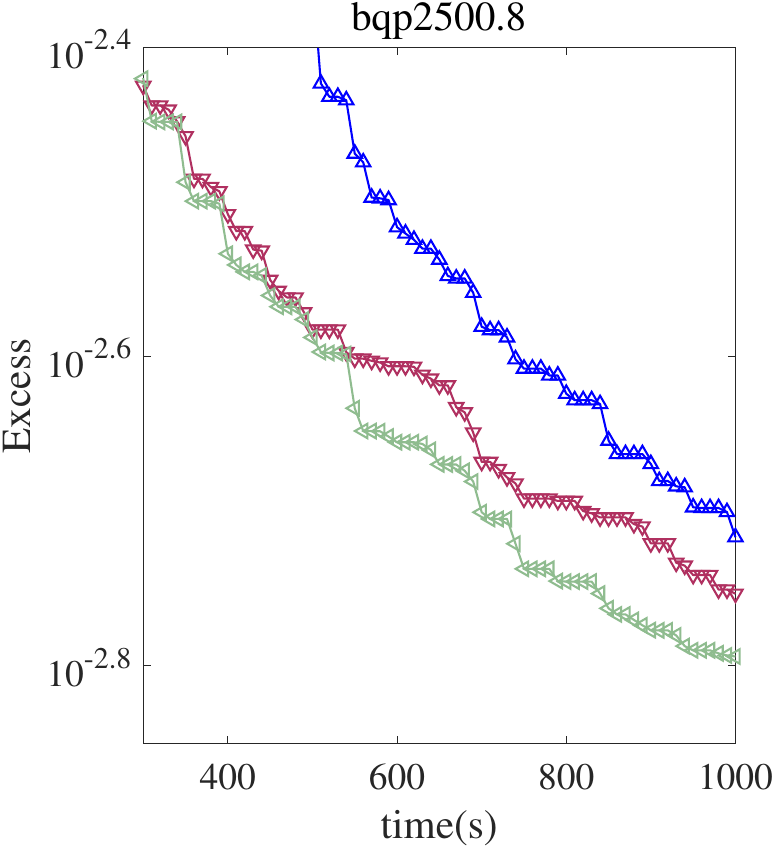}
		\label{se.bqp2500.8}}
	
	\subfloat[]{\includegraphics[width=1.8in]{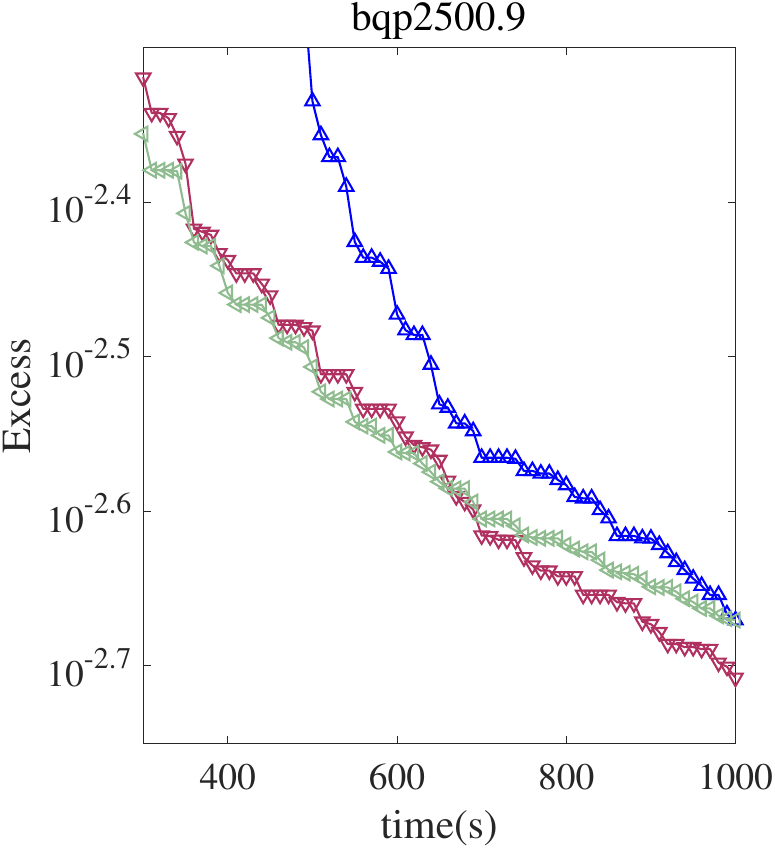}
		\label{se.bqp2500.9}}
	\subfloat[]{\includegraphics[width=1.8in]{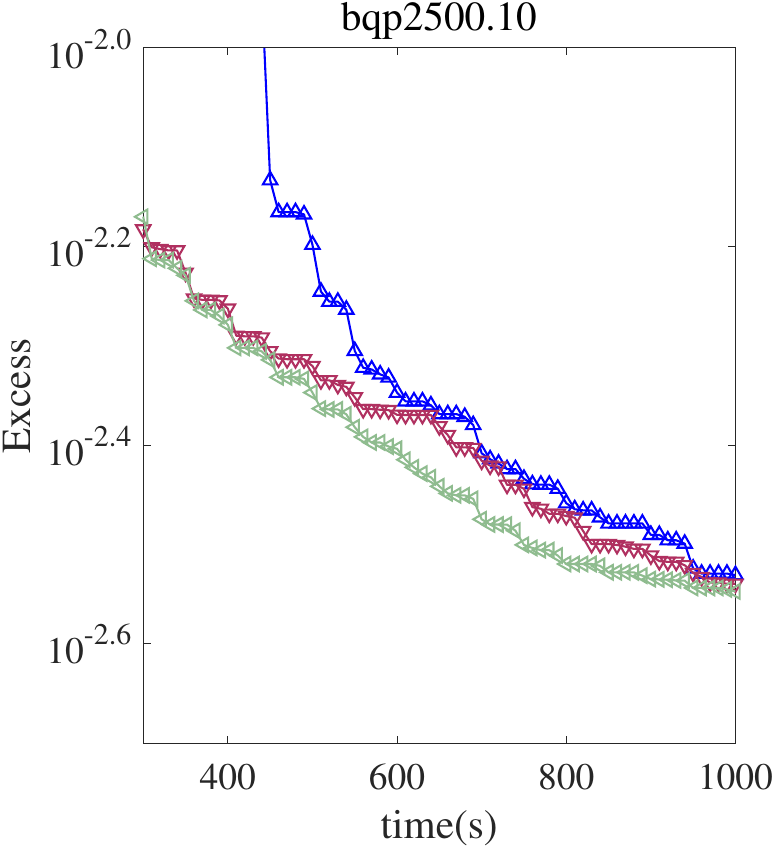}
		\label{se.bqp2500.10}}
	\caption{Comparison results of ILS, LSILS. (a) bqp2500.1. (b) bqp2500.2. (c) bqp2500.3. (d) bqp2500.4. (e) bqp2500.5. (f) bqp2500.6. (g) bqp2500.7. (h) bqp2500.8. (i) bqp2500.9. (j) bqp2500.10.}
	\label{seUBQP}
\end{figure*}

The result of the experiment is presented in \figurename~\ref{seUBQP}. The result shows that in all instances except bqp2500.5 and bqp2500.9, LSILS performs better than ILS throughout the whole search, \rv{which means that the HC transformation method reduces the difficulty of solving them.} For instance bqp2500.9, LSILS has a better effect than ILS before time point 600s, indicating that LSILS is still a promising landscape smoothing algorithm. LSILS has a lower curve than GH in all problems except bqp2500.1 and bqp2500.5, it implies that GH landscape smoothing method, first proposed for the TSP, is not applicable to the UBQP as HC transformation. In bqp2500.1 and bqp2500.5, LSILS performs better than GH in more than half the search time. Also, in 8 of 10 problems (except bqp2500.1 and bqp2500.6), GH has a higher excess curve than ILS for the whole process, reflecting that landscape smoothing effect of GH is not only worse than HC transformation, but even has a negative effect on metaheuristic. From the definition of GH transformation, the smoothing factor $\alpha$ is not equal to 1 in the first six local search rounds. When $\alpha$ is an even number, all elements of the smoothed $\widetilde{Q}$ is greater than 0, this operation will have a great impact on the landscape of the original UBQP, resulting in a negative effect on the landscape smoothing of the UBQP. From this result, it can be found that HC transformation is effective in smoothing the UBQP landscape.

\subsection{Performance of Parallel Cooperative LSILS for UBQPs}
According to previous analyses, HC transformation can smooth the landscapes of both UBQPs and TSPs, making these two NP-hard problems easier to solve. To use multi-core computers' capability to increase the performance of LSILS algorithm and further improve the smoothing effect of HC transformation, we conduct experiments on both UBQPs and TSPs using PC-LSILS that consists of multiple sequential LSILS algorithms running on the available cores. Additionally, the proposed PC-LSILS has a cooperation mechanism following a torus topology. To show the effectiveness of the PC-LSILS, we also introduce PI-LSILS which does not have cooperation among processes. Instances used in the experiments are the same as the aforementioned instances. About the compared algorithms, PC-LSILS is compared with PI-ILS (parallel version of ILS without cooperation and HC transformation), PI-LSILS (parallel version of LSILS without cooperation) and PI-GH (parallel version of GH without cooperation).

First, we tested PC-LSILS on UBQPs, where we set $m = 16$ processes and each process starts from a randomly and independently generated solution. At each time point (every 10 seconds for UBQPs), the minimum excess value among $m$ processes is selected as the excess of the parallel algorithm. The other experimental settings are the same as \ref{sectionLSILS}.

\begin{figure*}[!t]
	\subfloat[]{\includegraphics[width=1.8in]{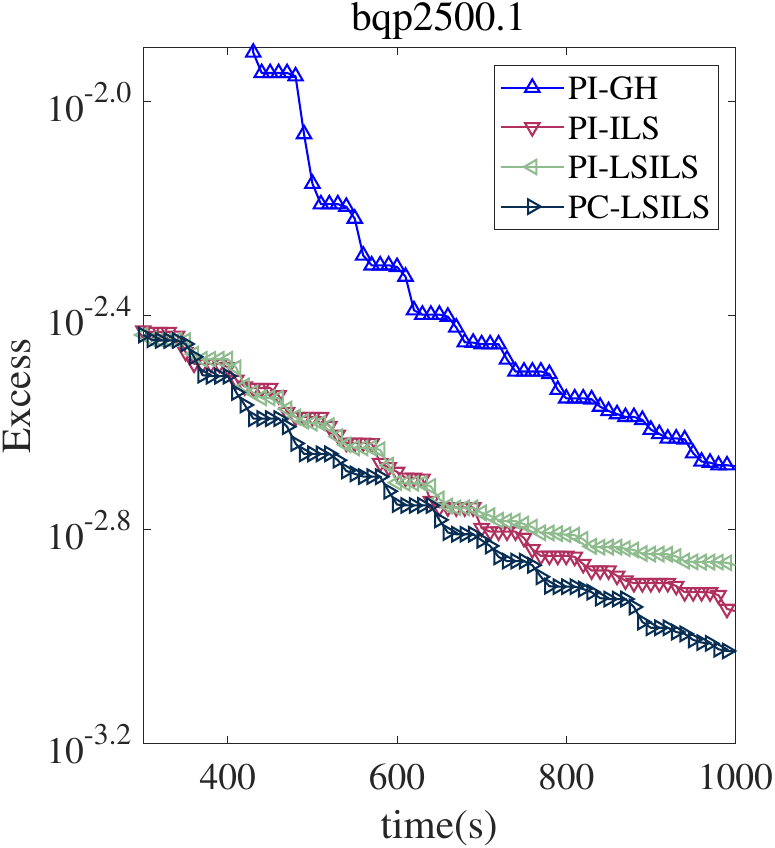}
		\label{bqp2500.1}}
	\subfloat[]{\includegraphics[width=1.8in]{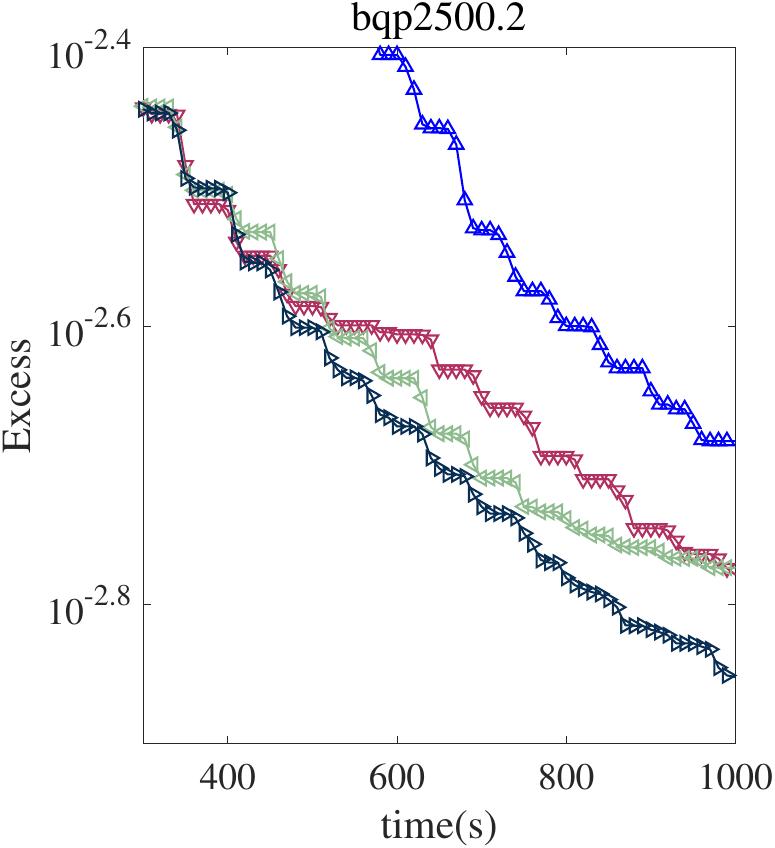}
		\label{bqp2500.2}}
	\subfloat[]{\includegraphics[width=1.8in]{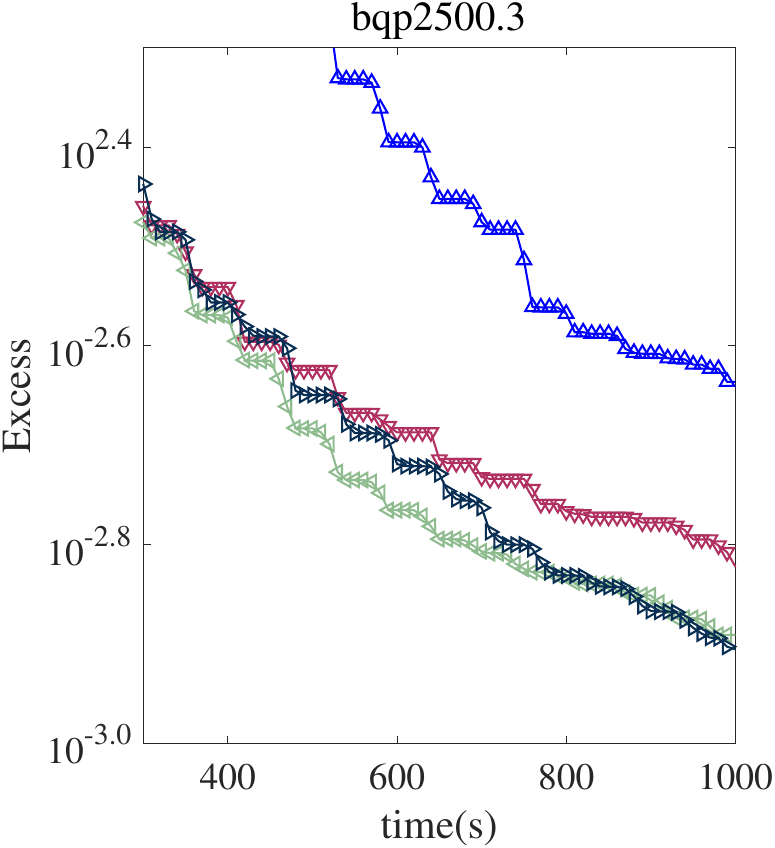}
		\label{bqp2500.3}}
	\subfloat[]{\includegraphics[width=1.8in]{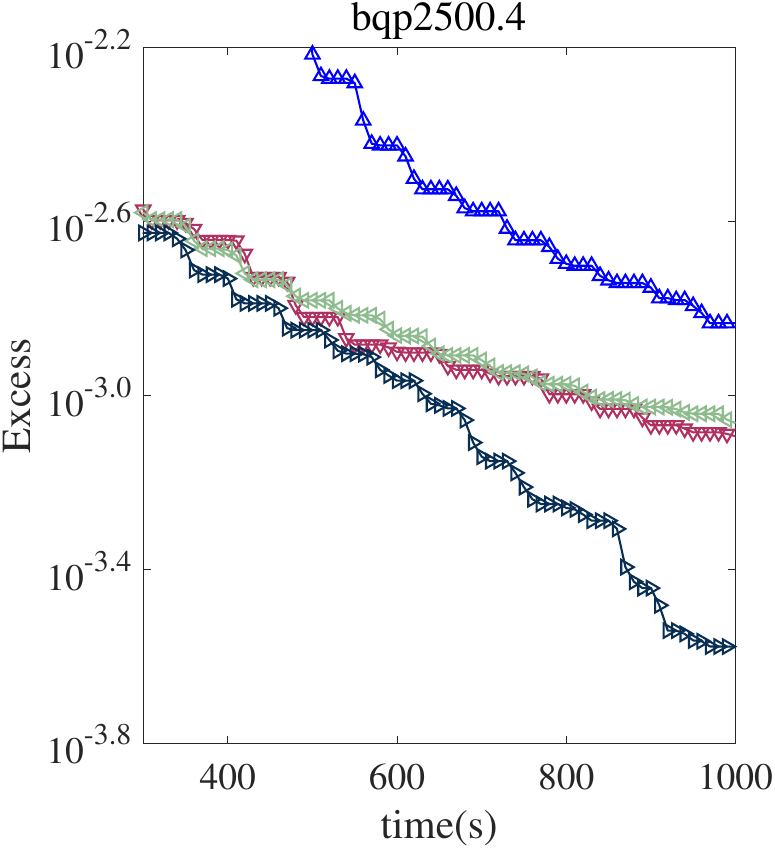}
		\label{bqp2500.4}}
	
	\subfloat[]{\includegraphics[width=1.8in]{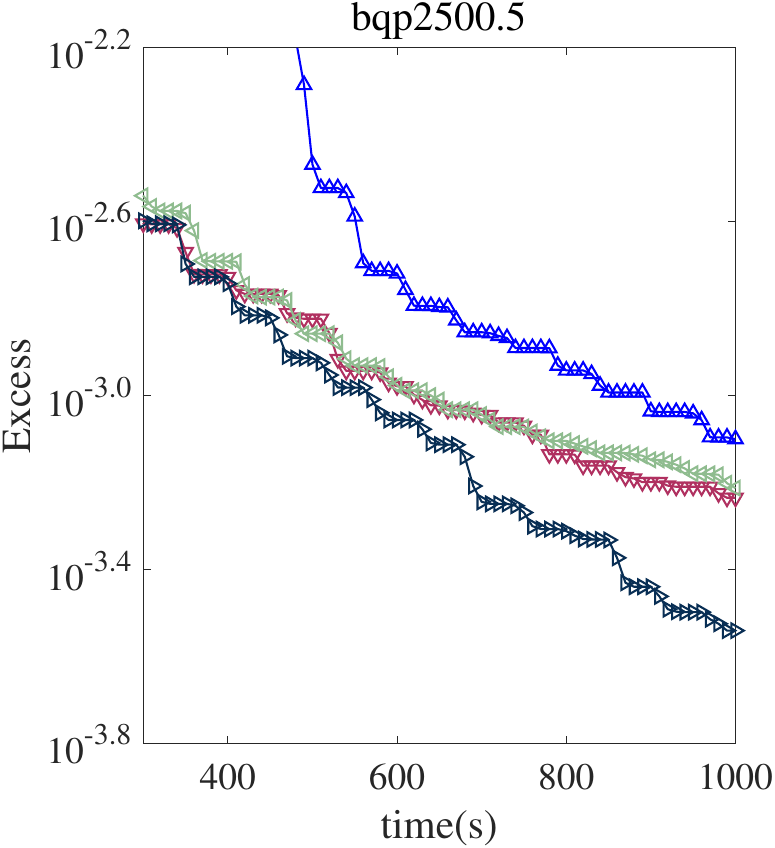}
		\label{bqp2500.5}}
	\subfloat[]{\includegraphics[width=1.8in]{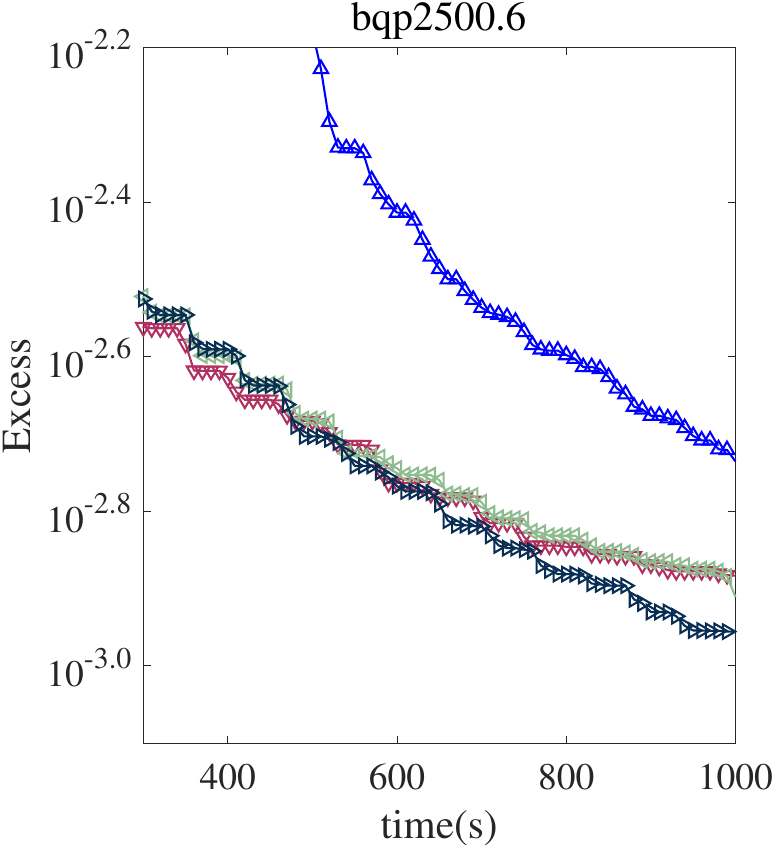}
		\label{bqp2500.6}}
	\subfloat[]{\includegraphics[width=1.8in]{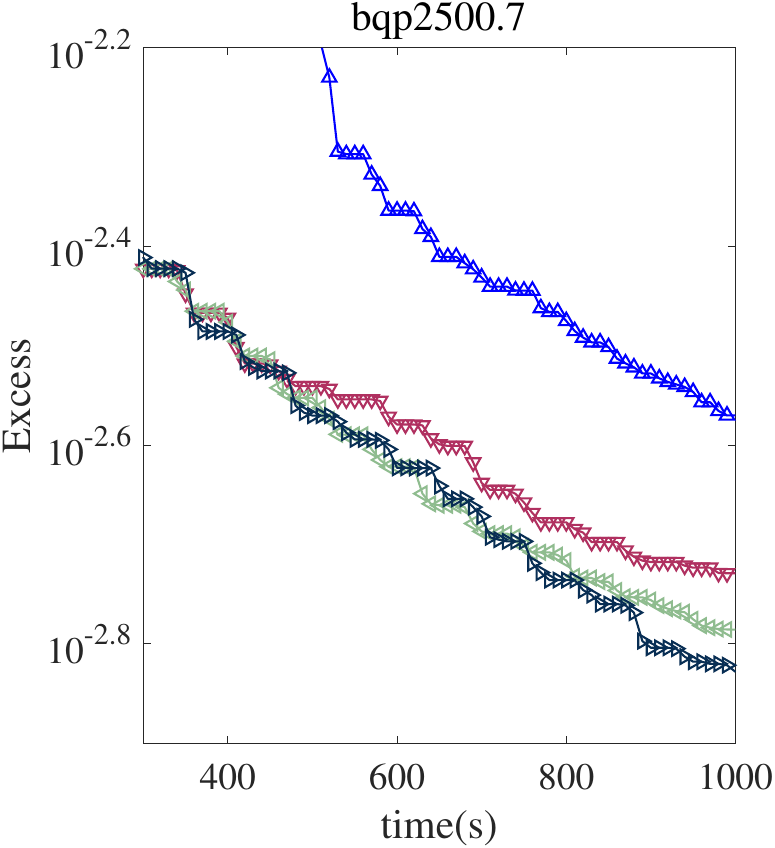}
		\label{bqp2500.7}}
	\subfloat[]{\includegraphics[width=1.8in]{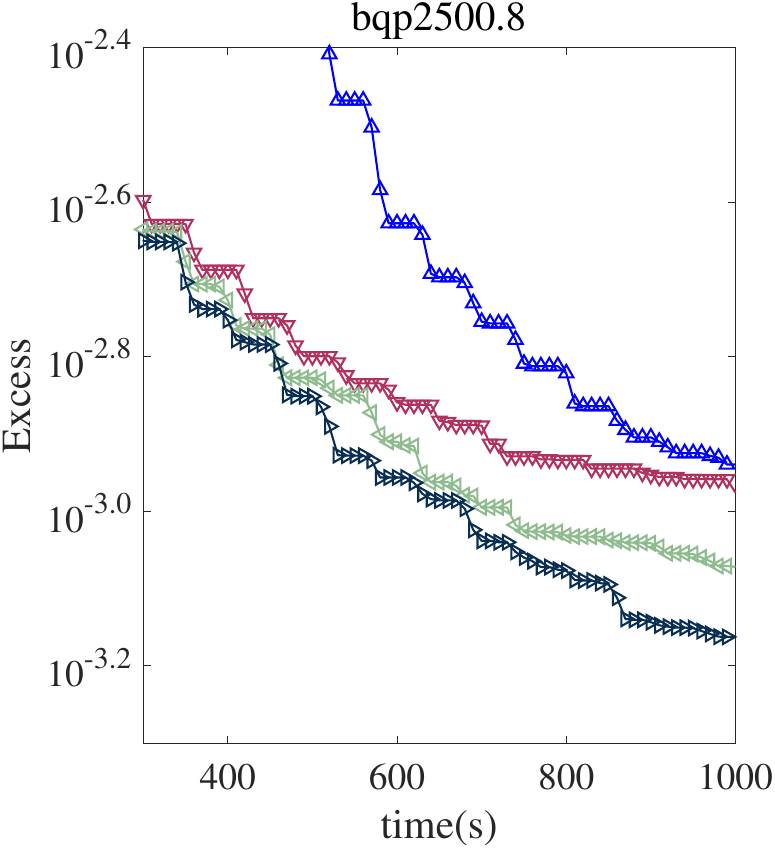}
		\label{bqp2500.8}}
	
	\subfloat[]{\includegraphics[width=1.8in]{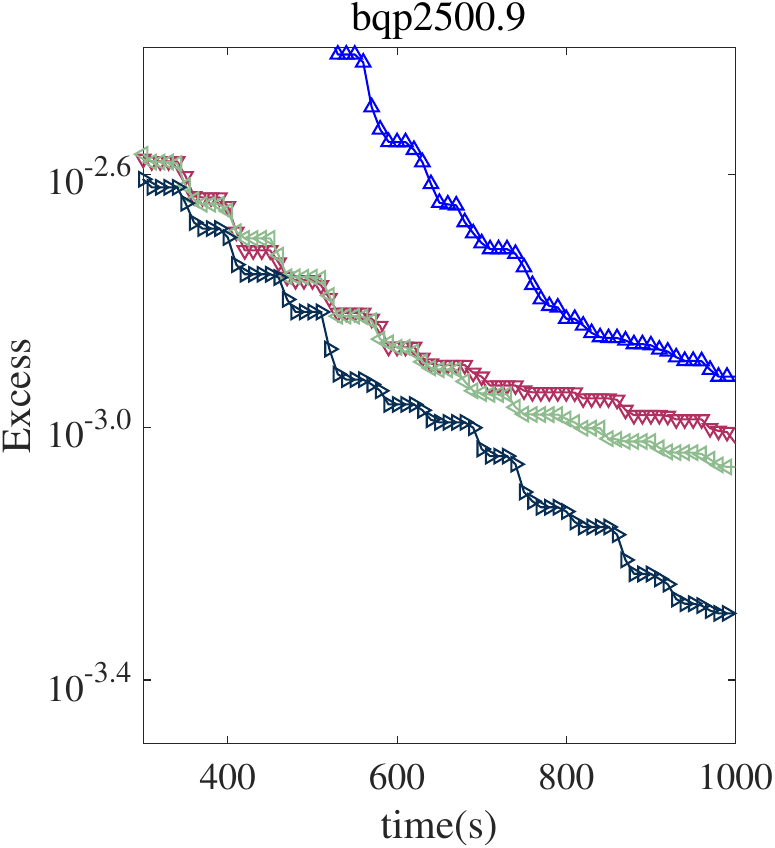}
		\label{bqp2500.9}}
	\subfloat[]{\includegraphics[width=1.8in]{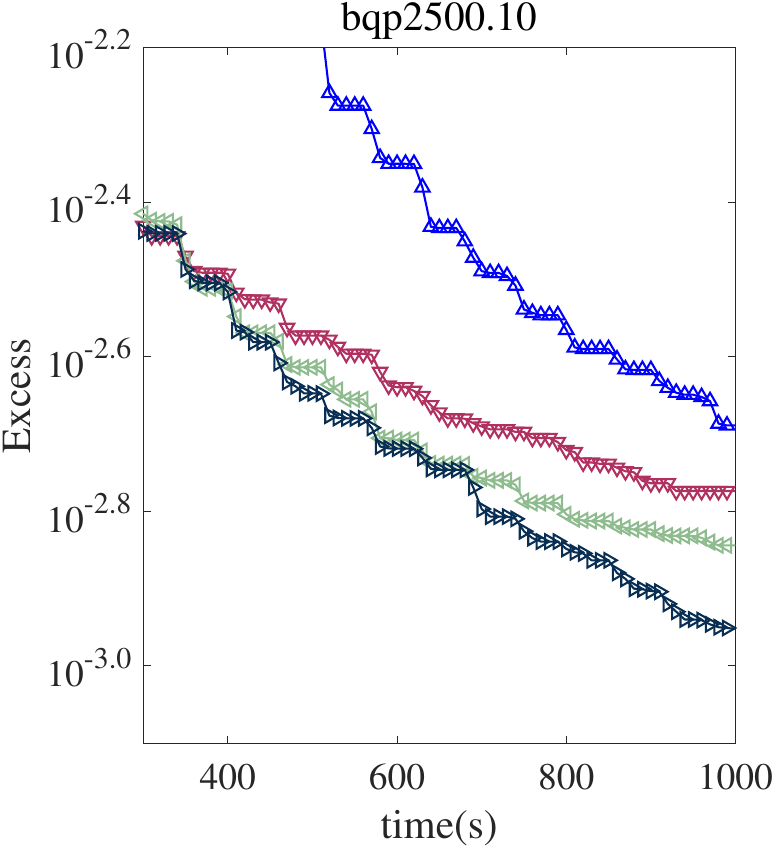}
		\label{bqp2500.10}}
	\caption{Comparison results of PI-GH, PI-ILS, PI-LSILS and PC-LSILS. (a) bqp2500.1. (b) bqp2500.2. (c) bqp2500.3. (d) bqp2500.4. (e) bqp2500.5. (f) bqp2500.6. (g) bqp2500.7. (h) bqp2500.8. (i) bqp2500.9. (j) bqp2500.10.}
	\label{parallelUBQP}
\end{figure*}

The result of this experiment is reported in \figurename~\ref{parallelUBQP}, which shows the average excess values of the best-so-far solutions found by different algorithms changed over time on. We can see that PI-LSILS outperforms both PI-ILS and PI-GH in 7 out of 10 instances with lower excess. This result indicates that parallel ILS with HC transformation can truly improve the performance of parallel ILS without HC transformation. It can be seen that these curves except PI-GH almost overlap before 400s and there is not much difference in effect of the algorithms, while after the time point 400s, the curves gradually separate. It is because at the beginning of the search, the quality of local optimal solution is low, and as the search and the exchange of elite solutions continue, high-quality local optimum is found, which improves the smoothing effect of HC transformation. \figurename~\ref{parallelUBQP} also shows that algorithm with HC transformation will have a more positive effect than with GH transformation.

Note that parallel LSILS with cooperation (PC-LSILS) has the lowest excess curve in all 10 instances, which reveals that cooperation among processes can indeed improve the search procedure.

\subsection{Performance of Parallel Cooperative LSILS for TSPs}\label{PC_LSILS_TSP}
In this section, we execute experiment on TSPs to illustrate that PC-LSILS can improve the smoothing effect of HC transformation. \rv{The HC transformation is first applied to TSPs by \cite{shi2020homotopic}. For a TSP with $n$ cities, they first construct a convex-hull TSP in which the $n$ cities lie on the convex hull of the TSP graph and the order of the cities follows the order they appear in a known local optimum of the original TSP. The convex-hull TSP can be proved unimodal. Then the original TSP can then be smoothed by a convex combination of the convex-hull TSP and the original TSP with a parameter $\lambda \in$ [0,1].}

Different from UBQP experiments, in this experiment, the stopping criteria is set 60 seconds of CPU runtime. We test the performance of the proposed parallel framework based on 3-Opt local search \citep{cook2011traveling} with first-improvement strategy and double bridge perturbation. The developed PC-LSILS and PI-LSILS are compared against PI-ILS, PI-GH and PI-SSA.  \rv{The compared landscape smoothing algorithm GH is proposed by \cite{gu1994efficient} and SSA is proposed by \cite{coy1998see}. We have introduced the landscape smoothing process of GH in equation (\ref{dij}). For SSA, first, the distances between cities are smoothed with a convex function and the local search heuristic is applied. Second, the distances are smoothed with a concave function and the local search heuristic is applied. The convex function and the concave function are shown as Eq.~ (\ref{SSA}).
\begin{equation}
\label{SSA}
\begin{aligned}
&\textrm{Concave}:
d_{ij}(\alpha) = d_{ij}^{1/\alpha},\\
&\textrm{Convex}:
d_{ij}(\alpha) = d_{ij}^{\alpha}.
\end{aligned}
\end{equation}
In the experiment, for GH, a sequence of $\alpha$ values ranging from 6 to 1 is employed across the initial six local search rounds. Similarly, for SSA, a sequence comprising $\alpha$ values of 7, 5, 3, and 1 is utilized during the first four local search rounds. In both GH and SSA, once $\alpha$ reaches 1 and the function evaluation budget has not been run out, the algorithm will apply ILS on the original TSP landscape, characterized by $\alpha = 1$, until the completion of the function evaluation budget.}

We parallelize GH and SSA and then compare them with our parallel LSILS. For PI-LSILS and PC-LSILS, we set two different settings of $\lambda$, with constant $\lambda$ or dynamically changing $\lambda$. For dynamic $\lambda$, it increases from 0 at the beginning to 0.09 during the search. For example, parameter $\lambda$ will change from 0 to 0.01 at 6 s, from 0.01 to 0.02 at 12 s, and so on. The other experimental settings are the same as UBQPs.

\begin{figure*}[!t]
	\subfloat[]{\includegraphics[width=1.8in]{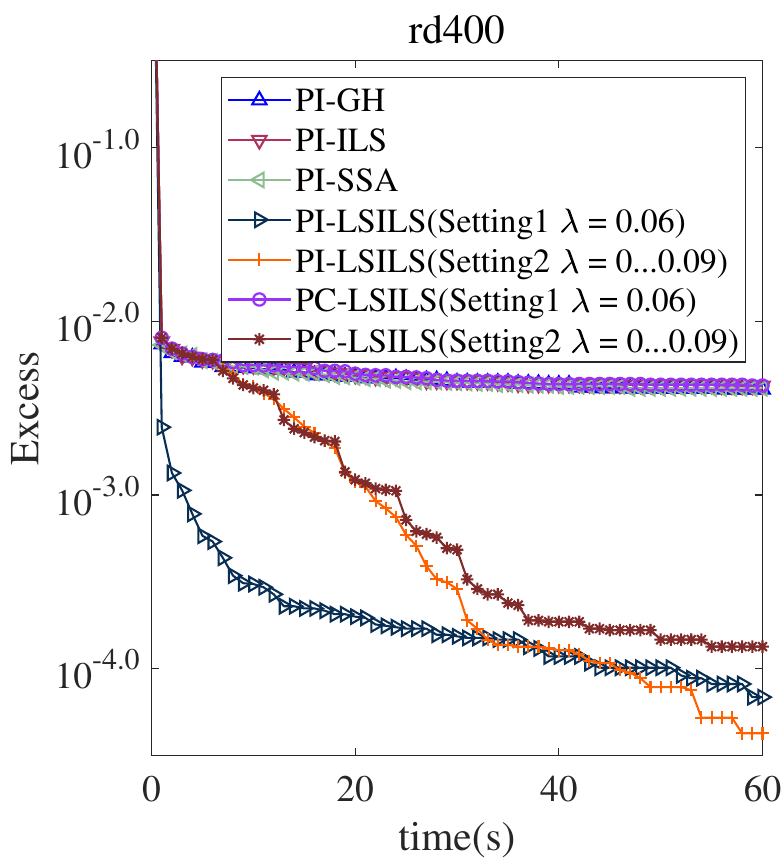}
		\label{rd400}}
	\subfloat[]{\includegraphics[width=1.8in]{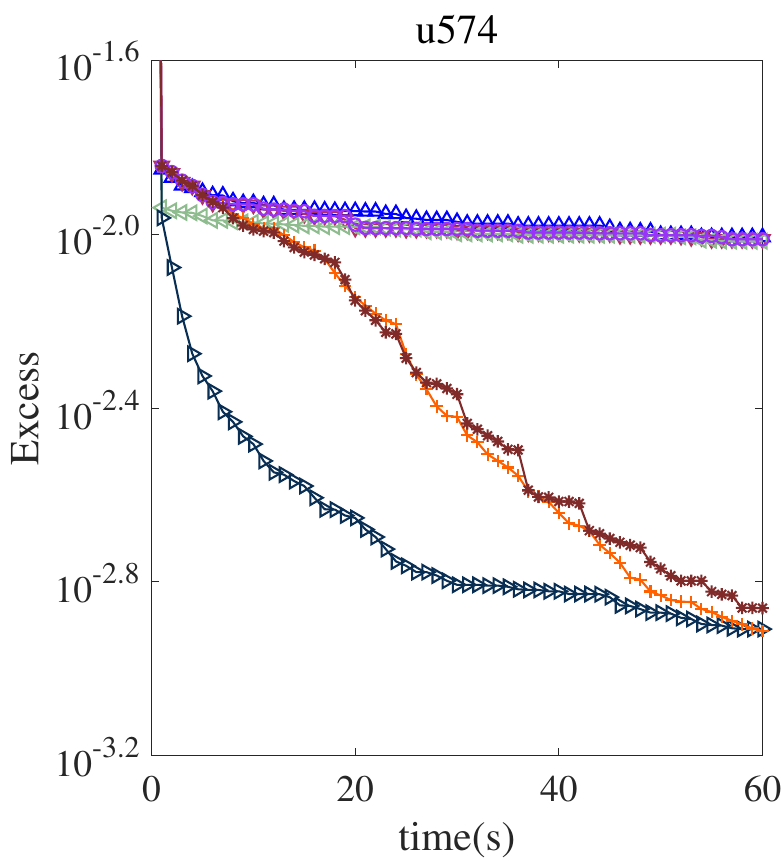}
		\label{u574}}
	\subfloat[]{\includegraphics[width=1.8in]{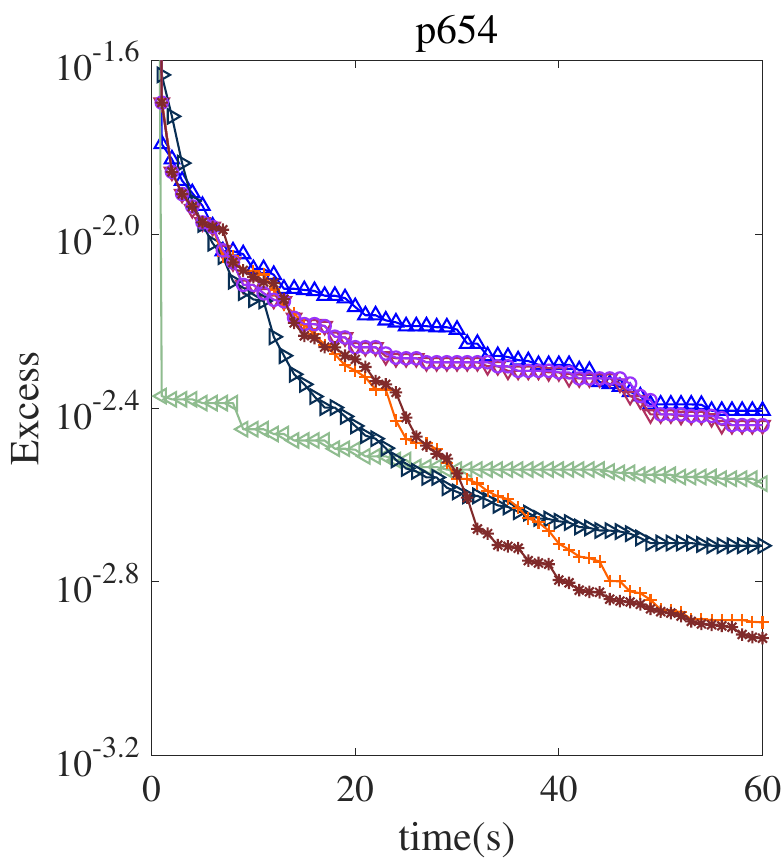}
		\label{p654}}
	\subfloat[]{\includegraphics[width=1.8in]{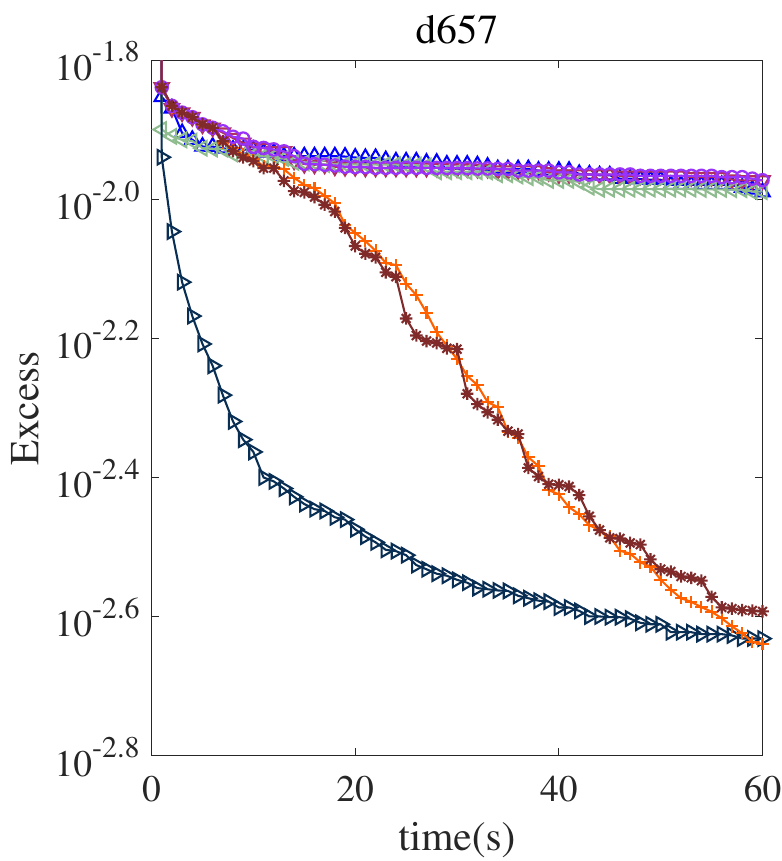}
		\label{d657}}
	
	\subfloat[]{\includegraphics[width=1.8in]{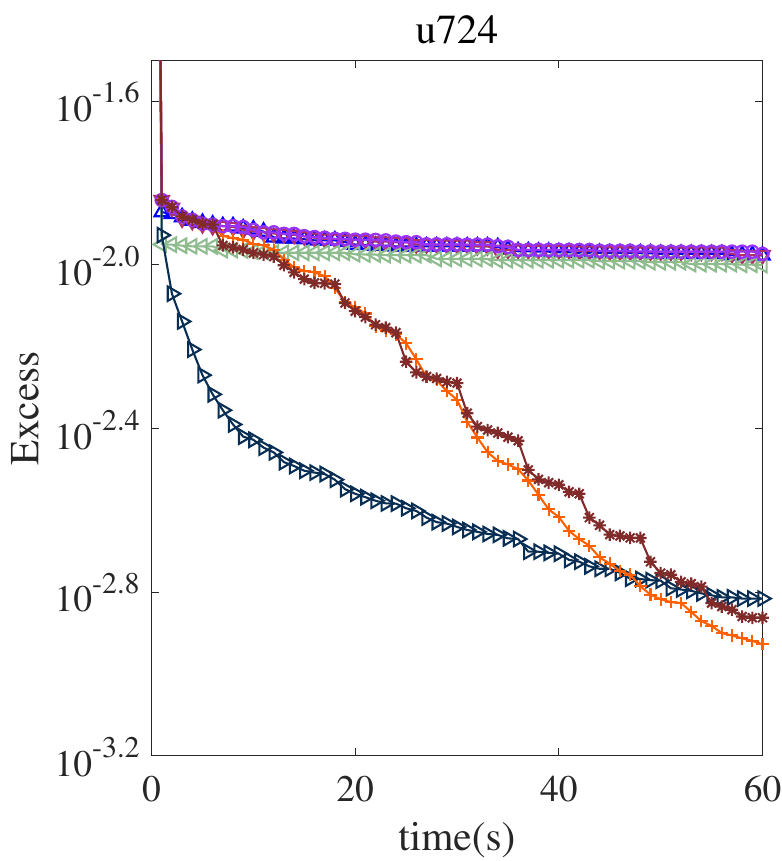}
		\label{u724}}
	\subfloat[]{\includegraphics[width=1.8in]{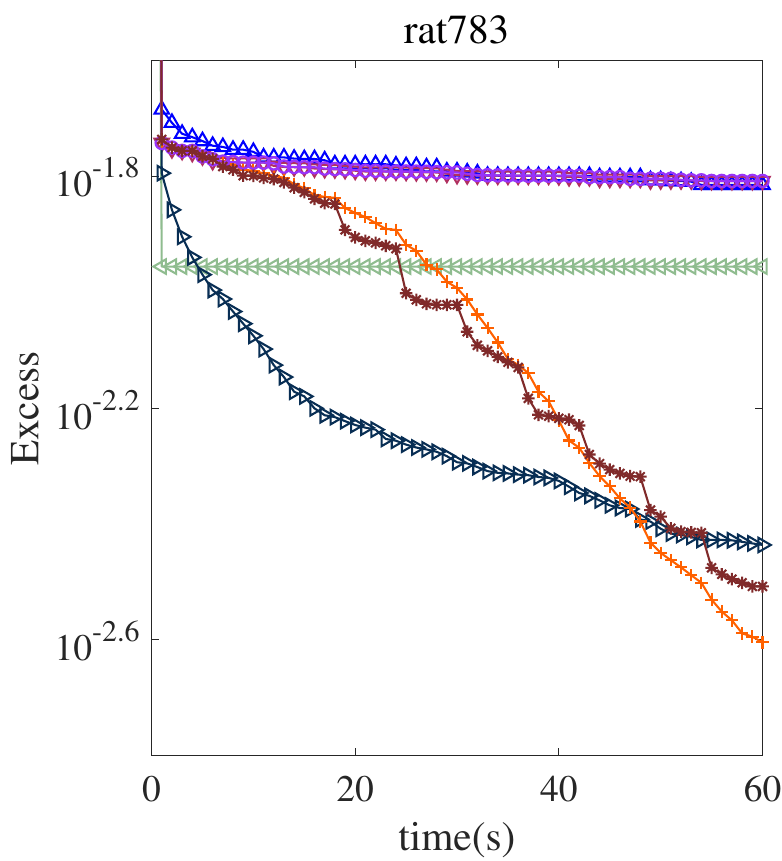}
		\label{rat783}}
	\subfloat[]{\includegraphics[width=1.8in]{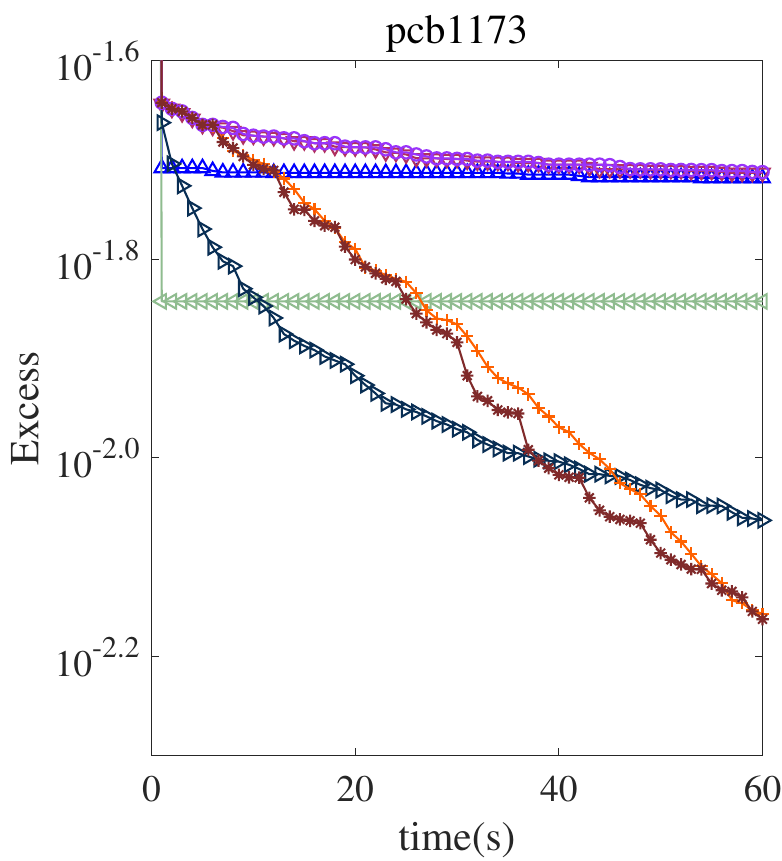}
		\label{pcb1173}}
	\subfloat[]{\includegraphics[width=1.8in]{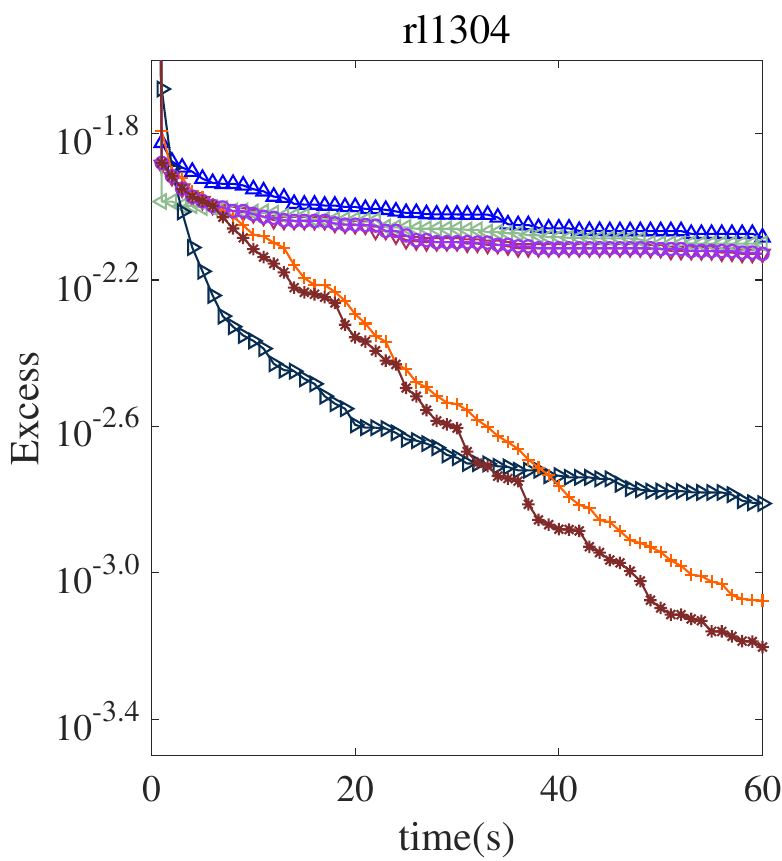}
		\label{rl1304}}
	
	\subfloat[]{\includegraphics[width=1.8in]{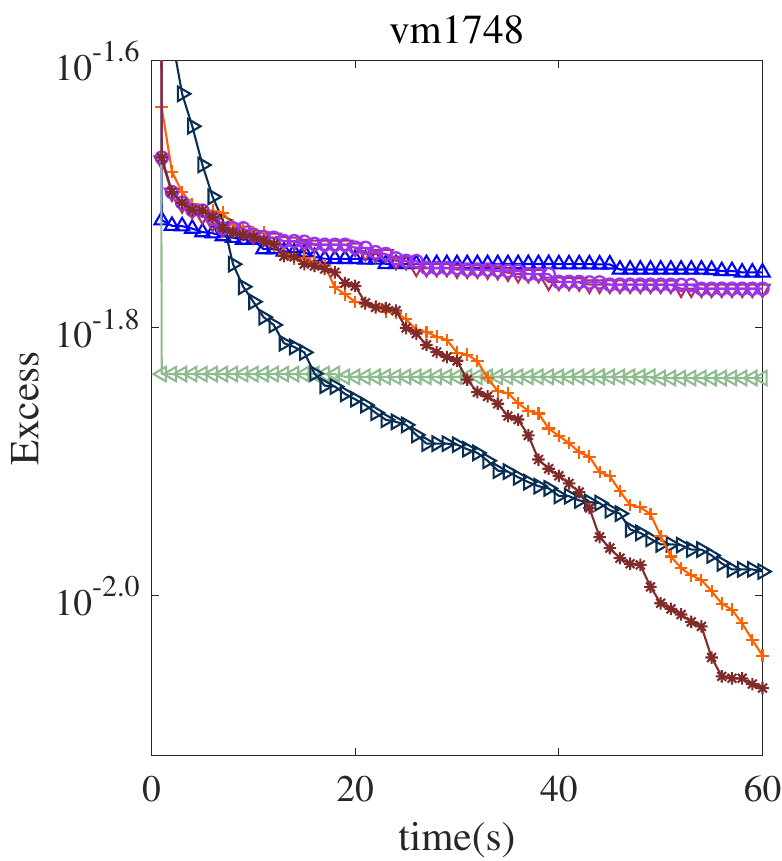}
		\label{vm1748}}
	\subfloat[]{\includegraphics[width=1.8in]{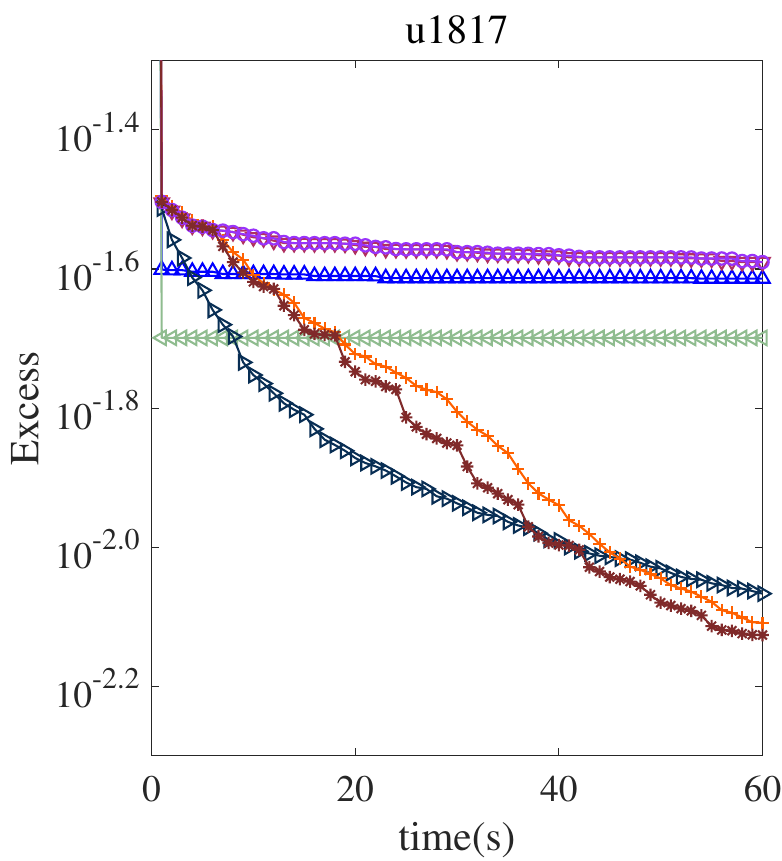}
		\label{u1817}}
	\caption{Comparison results of PI-GH, PI-ILS, PI-SSA, PI-LSILS and PC-LSILS based on 3-Opt local search and double bridge perturbation. (a) rd400. (b) u574. (c) p654. (d) d657. (e) usa13509. (f) u724. (g) rat783. (h) rl1304. (i) vm1748. (j) u1817.}
	\label{parallelTSP}
\end{figure*}

\figurename~\ref{parallelTSP} shows the comparison results with other algorithms, which demonstrates the average excess values over time. It can be observed that in 7 out of 10 test instances, PI-LSILS with setting2 and PC-LSILS with setting2 get the two lowest excess value among these parallel algorithms which means that parallel LSILS with increasing $\lambda$ values achieves the best performance. Meanwhile, three parallel algorithms (PI-LSILS with setting1, PI-LSILS with setting2, PC-LSILS with setting2) perform significantly better than PI-ILS, PI-GH, and PI-SSA in terms of excess. This result clearly indicates that the parallel LSILS algorithm can improve the efficiency for solving TSPs. Additionally, at the beginning of the search, the performance of parallel LSILS is also relatively poor. As explained before, this is because that the smoothing effect of the HC transformation depends highly on the quality of the local optimum used to construct the convex-hull TSP. As the search goes by, better solutions are found, the performance of parallel LSILS becomes better and exceeds the performance of PI-GH, PI-ILS and PI-SSA eventually.

Actually, in relatively small instances (rd400 , u574, d657, u724, rat783), the minimum excess value in the end is obtained by PI-LSILS with setting2. Conversely, in relatively large-size instances, specifically p654, pcb1173, rl1304, vm1748, and u1817, PC-LSILS with setting2 performs best. This result indicates that the technique of communication among processes has a better effect in larger-size TSPs. A possible explanation for this phenomenon is that local search can easily find global optimal solutions for smaller TSPs and thus sharing gathered elite solutions among the processes will be less effective than only executing local search. However, in larger-size instances where the solution space is much more complex, so it is not so efficient to conduct just local search. Through sharing the information of elite solutions among processes, a process can get better solutions to construct the unimodal UBQP, resulting in a better algorithm experimental results.

However, the excess curve of PC-LSILS with setting1 is higher than PI-LSILS, even higher than both PI-SSA and PI-GH, meaning that cooperation in PC-LSILS with setting1 has negative influence while has positive influence with setting2. We assume that it is because at the beginning of search, the local optima have relatively low quality, but PC-LSILS spreads them to other processes, resulting in all processes being attracted by low-quality solutions, which has poor HC smoothing effect, thereby reducing the performance of the parallel algorithm. But with PI-LSILS (setting1), each process conducts its own search, which increases the diversity of the search and avoids performing HC transformations on poor-quality local optimal solutions.

In summary, the proposed parallel cooperative ILS with HC transformation based on local optima can indeed improve the global search ability of  algorithms on both UBQPs and TSPs. Elite solutions sharing seems to perform better than the independent search approach and is a valuable strategy to tackle large-size instances, while the advantage is not obvious in small-size instances.

\section{Conclusion}\label{6}
In this paper, a landscape smoothing method HC transformation is applied to UBQPs, which is defined as a homotopic combination of the original UBQP and a toy UBQP. Importantly, this paper provides a rigorous proof confirming that the toy UBQP is unimodal. \rv{To investigate the smoothing effects of the HC transformation based on global or local optima, we conduct landscape analysis experiments with two metrics to characterize the ruggedness. The experimental results show that the HC transformation based on global or local optima can both smooth the UBQP landscape.} Then LSILS, a landscape smoothing-based iterative algorithm using HC transformation is proposed for UBQPs. Furthermore, based on LSILS, we introduce a parallel cooperative algorithm framework (PC-LSILS) for both TSPs and UBQPs. Experimental results show the potential of HC transformation in addressing UBQPs with rugged fitness landscapes and highlight the superior performance of the parallel cooperative LSILS approach.

Looking forward, our future research endeavors will involve applying the proposed HC transformation technique to the Quadratic Assignment Problem (QAP). A key challenge in this context will be the design of a unimodal QAP based on local optima, and we will explore innovative approaches to address this challenge.

\bibliographystyle{model5-names}
\bibliography{UBQP_parallel_smoothing}

\end{document}